\documentclass[11pt]{amsart}
\usepackage{amsfonts}
\usepackage{amsmath, amsthm, amssymb,color}
\usepackage{latexsym}
\newcommand{\vect}[1]{\boldsymbol{\mathbf{#1}}}
\numberwithin{equation}{section}
\allowdisplaybreaks

\makeatletter
\renewcommand\@biblabel[1]{}
\makeatother
\makeatletter
\CheckCommand*{\@cite}[2]{%
  {%
    \@citestyle[\citeform{#1}\if@tempswa, #2\fi]%
  }%
}
\renewcommand*{\@cite}[2]{%
  {%
    \@citestyle\citeform{#1}\if@tempswa, #2\fi
  }%
}
\makeatother

\usepackage{graphics}
\usepackage{graphicx}
\usepackage{epsfig}
\usepackage{mathrsfs}

\definecolor{darkblue}{rgb}{.1, 0.1,.8}
\definecolor{darkgreen}{rgb}{0,0.8,0.2}
\definecolor{darkred}{rgb}{.8, .1,.1}
\newcommand{\E}{{\mathbb E}}
\renewcommand{\P}{{\mathbb P}}
\newcommand{\bfi}{\begin{fig}}
\newcommand{\efi}{\end{fig}}
\newcommand{\levy}{L\'evy}

\newcommand{\clt}{central limit theorem}

\newcommand{\ex}{{\rm e}\,}

\newcommand{\asy}{asymptotic}

\newcommand{\id }{infinitely divisible}

\newcommand{\ts}{time series}
\newcommand{\tsa}{\ts\ analysis}

\newtheorem{lemma}{Lemma}[section]

\newtheorem{theorem}[lemma]{Theorem}

\newtheorem{proposition}[lemma]{Proposition}
\newtheorem{definition}[lemma]{Definition}
\newtheorem{corollary}[lemma]{Corollary}
\newtheorem{example}[lemma]{Example}
\newtheorem{exercise}[lemma]{Exercise}
\newtheorem{remark}[lemma]{Remark}
\newtheorem{fig}[lemma]{Figure}
\newtheorem{tab}[lemma]{Table}

\newcommand{\cid}{\stackrel{d}{\rightarrow}}
\newcommand{\cip}{\stackrel{p}{\rightarrow}}

\newcommand{\bfQ}{{\bf Q}}

\newcommand{\bth}{\begin{theorem}}
\newcommand{\ethe}{\end{theorem}}

\newcommand{\bre}{\begin{remark}\em }
\newcommand{\ere}{\end{remark}}

\newcommand{\ble}{\begin{lemma}}
\newcommand{\ele}{\end{lemma}}

\newcommand{\bde}{\begin{definition}}
\newcommand{\ede}{\end{definition}}
\newcommand{\bco}{\begin{corollary}}
\newcommand{\eco}{\end{corollary}}

\newcommand{\bpr}{\begin{proposition}}
\newcommand{\epr}{\end{proposition}}

\newcommand{\bpf}{\begin{proof}}
\newcommand{\epf}{\end{proof}}

\newcommand{\bexer}{\begin{exercise}}
\newcommand{\eexer}{\end{exercise}}

\newcommand{\bexam}{\begin{example}\rm }
\newcommand{\eexam}{\end{example}}

\newcommand{\btab}{\begin{tab}}
\newcommand{\etab}{\end{tab}}

\newcommand{\lhs}{left-hand side}
\newcommand{\fidi}{finite-dimensional distribution}
\newcommand{\rv}{random variable}

\newcommand{\var}{{\rm Var}}

\newcommand{\cov}{{\rm cov}}

\newcommand{\as}{{\rm a.s.}}

\newcommand{\rhs}{right-hand side}

\newcommand{\beao}{\begin{eqnarray*}}
\newcommand{\eeao}{\end{eqnarray*}\noindent}

\newcommand{\beam}{\begin{eqnarray}}
\newcommand{\eeam}{\end{eqnarray}\noindent}

\newcommand{\beqq}{\begin{equation}}
\newcommand{\eeqq}{\end{equation}\noindent}

\newcommand{\beqo}{\begin{equation*}}
\newcommand{\eeqo}{\end{equation*}\noindent}

\newcommand{\bce}{\begin{center}}
\newcommand{\ece}{\end{center}}

\newcommand{\barr}{\begin{array}}
\newcommand{\earr}{\end{array}}

\newcommand{\stp}{\stackrel{\P}{\rightarrow}}
\newcommand{\std}{\stackrel{d}{\rightarrow}}
\newcommand{\stas}{\stackrel{\rm a.s.}{\rightarrow}}

\newcommand{\eqd}{\stackrel{d}{=}}

\newcommand{\vague}{\stackrel{\lower0.2ex\hbox{$\scriptscriptstyle
                    \it{v} $}}{\rightarrow}}
\newcommand{\weak}{\stackrel{\lower0.2ex\hbox{$\scriptscriptstyle
                    \it{w} $}}{\rightarrow}}
\newcommand{\what}{\stackrel{\lower0.2ex\hbox{$\scriptscriptstyle
                    \it{\hat{w}} $}}{\rightarrow}}

\newcommand{\bdis}{\begin{displaymath}}
\newcommand{\edis}{\end{displaymath}\noindent}

\newcommand{\R}{\mathbb{R}}
\newcommand{\Xbold}{{\mathbf{X}}}

\newcommand{\nto}{n\to\infty}

\newcommand{\ov}{\overline}
\newcommand{\wt}{\widetilde}
\newcommand{\wh}{\widehat}
\newcommand{\vep}{\varepsilon}

\newcommand{\regvary}{regularly varying}

\newcommand{\bbr}{{\mathbb R}}

\newcommand{\bbz}{{\mathbb Z}}

\newcommand{\bbs}{{\mathbb S}}

\newcommand{\con}{convergence}

\newcommand{\st}{such that}
\newcommand{\fif}{if and only if}
\newcommand{\wrt}{with respect to}
\newcommand{\chf}{characteristic function}
\newcommand{\fct}{function}

\newcommand{\ds}{distribution}

\newcommand{\rep}{representation}

\newcommand{\seq}{sequence}

\newcommand{\pro}{probabilit}

\newcommand{\ms}{measure}

\newcommand{\bfX}{{\bf X}}

\renewcommand\d{{\mathrm d}}

\newcommand\1{{\mathbf 1}}

\newcommand{\bealm}{\begin{align}}
\newcommand{\eealm}{\end{align}\noindent}

\begin{document}

\bibliographystyle{alpha}
\title[Applications of Distance Correlation to Time Series]{Applications of Distance Correlation to Time Series}
\today
\author[R.A. Davis]{Richard A. Davis}
\address{Department of Statistics, Columbia University, 1255 Amsterdam Ave, New York, NY 10027, USA.}
\email{rdavis@stat.columbia.edu}
\author[M. Matsui]{Muneya Matsui}
\address{Department of Business Administration, Nanzan University, 18
Yamazato-cho, Showa-ku, Nagoya 466-8673, Japan.}
\email{mmuneya@gmail.com}
\author[T. Mikosch]{Thomas Mikosch}
\address{Department of Mathematics,
University of Copenhagen,
Universitetsparken~5,
DK-2100 Copenhagen,
Denmark}
\email{mikosch@math.ku.dk}
\author[P. Wan]{Phyllis Wan}
\address{Department of Statistics, Columbia University, 1255 Amsterdam Ave, New York, NY 10027, USA.}
\email{phyllis@stat.columbia.edu}
\begin{abstract}
 The use of empirical characteristic functions for inference problems,
 including estimation in some special parametric settings and testing
 for goodness of fit, has a long history dating back to the 70s (see for
 example, \cite{feuerverger:mureika:1977}, Cs\"{o}rg\H{o}\ (1981a,1981b,1981c)
 \cite{feuerverger:1993}). More recently, there has been
 renewed interest in using empirical characteristic functions in other
 inference settings. The distance covariance and correlation, developed
 by \cite{szekely:rizzo:2009} for measuring
 dependence and testing independence between two random vectors, are
 perhaps the best known illustrations of this. We apply these ideas to
 stationary univariate and multivariate time series to measure lagged
 auto- and cross-dependence in a time series. Assuming strong mixing, we
 establish the relevant asymptotic theory for the sample auto- and
 cross-distance correlation functions.  We also apply the auto-distance
 correlation function (ADCF) to the residuals of an autoregressive
 processes as a test of goodness of fit.   Under the null that an
 autoregressive model is true, the limit distribution of the empirical
 ADCF
 can differ markedly from the corresponding one based on an iid
 sequence.  We illustrate the use of the empirical auto- and
 cross-distance correlation functions for testing dependence and
 cross-dependence of time series in a variety of different contexts. 
 
\end{abstract}
\keywords{Auto- and cross-distance correlation \fct , testing independence, \ts , strong mixing, ergodicity,
Fourier analysis, $U$-statistics, AR process, residuals}

\subjclass{Primary 62M10; Secondary 60E10 60F05 60G10 62H15 62G20}
\thanks{Thomas Mikosch's research is partly supported by the Danish Research Council Grant DFF-4002-00435 ``Large random matrices with heavy tails and dependence''. Richard Davis and Phyllis Wan's research is supported in part by ARO MURI grant W911NF-12-1-0385.  
Muneya Matsui's research is supported by the JSPS Grant-in-Aid
for Young Scientists B (16k16023).}
\maketitle
\section{Introduction}
In time series analysis,  modeling  serial dependence is typically
the overriding objective.  In order to achieve this goal, it is
necessary to formulate a measure of dependence and this may depend on
the features in the data that one is trying to capture.  The
autocorrelation function (ACF), which provides a measure of linear
dependence, is perhaps the most used dependence measure in time series.
It is closely linked with the class of ARMA models and provides guidance
in both model selection and model confirmation.  On the other hand, the
ACF gives only a partial description of serial dependence.  As seen with
financial time series, data are typically uncorrelated but dependent for
which the ACF is non-informative.  In this case, the dependence becomes
visible by examining the ACF applied to the absolute values or squares
of the time series.  In this paper, we consider the application of
distance correlation in a time series setting, which can overcome some
of the limitations with other dependence measures. 

In recent years, the notions of  distance covariance and correlation have become rather popular in applied statistics.  Given vectors $X$ and $Y$ with values in $\bbr^p$ and $\bbr^q$,
the {\em distance covariance} between $X$ and $Y$ \wrt\ a suitable \ms\ $\mu$ on $\bbr^{p+q}$ is given by
\beam\label{eq:i}
T(X,Y;\mu)= \int_{\bbr^{p+q}} \big|\varphi_{X,Y}(s,t)-\varphi_X(s)\,\varphi_Y(t)\big|^2\,\mu(ds,dt)\,,
\eeam
where we denote the \chf\ of any random vector $Z\in\bbr^d$ by
\beao
\varphi_Z(t) = \E [\ex^{i \langle t,Z\rangle}],\qquad t\in\bbr^d \,.
\eeao
The {\em distance correlation} is the corresponding version of $T$ standardized to values in $[0,1]$.
The quantity $T(X,Y;\mu)$ is zero \fif\ $\varphi_{X,Y}=\varphi_X\,\varphi_Y$ $\mu$-a.e. In many situations, for example when $\mu$
has a positive Lebesgue density on $\bbr^{p+q}$, we may conclude that $X$ and $Y$ are independent \fif\ $T(X,Y;\mu)=0$.
An empirical version $T_n(X,Y; \mu)$ of
$T(X,Y;\mu)$ is obtained if the \chf s in \eqref{eq:i} are replaced by their corresponding empirical
versions. Then one can build a test for independence between $X$ and $Y$
based on the \ds\ of $T_n$ under the null hypothesis that $X$ and $Y$ are independent.

The use of empirical \chf s for univariate and multivariate \seq s for
inference purposes has  a long history. In the 1970s and 1980s,
\cite{feuerverger:mureika:1977}, Cs\"{o}rg\H{o}\ (1981a,1981b,1981c)
and many others proved fundamental \asy\ results for iid \seq s,
including Donsker-type theory for the empirical \chf .
Statisticians have applied these methods for
goodness-of-fit tests, changepoint detection, testing for independence,
etc.; see for example Meintanis and coworkers
(\cite{meintanis:iliopoulos:2008,hlavka:kuskova:meintanis:2011,meintanis:ngatchou:taufer:2015}),
and the references therein.
The latter authors employed the empirical distance covariance for finite \ms s $\mu$. \cite{feuerverger:1993} was the first to apply statistics of form \eqref{eq:i} for general measures. In particular, he advocated the infinite measure
\beao\label{eq:feuermeas}
\mu(ds,dt)=|s|^{-2}|t|^{-2}ds\,dt\,
\eeao
for testing independence of univariate data.
Sz\'ekely et al.\footnote{They appeared to have coined the terms
distance covariance and correlation.}
(\cite{szekely:rizzo:bakirov:2007}, \\
\cite{szekely:rizzo:2014,szekely:rizzo:2009}, see also the references therein)
developed \asy\ techniques for the empirical distance covariance and correlation of
iid \seq s for the infinite measure  $\mu$ given by
\beam\label{eq:szikelymeas}
\mu(ds,dt)=c_{p,q}|s|^{-\alpha-p}|t|^{-\alpha-q}ds\,dt\,
\eeam
where $c_{p,q}$ is a constant (see \eqref{eq:density}) and $\alpha\in (0,2)$. With this choice of $\mu$, the distance correlation, $T(X,Y;\mu)/(T(X,X;\mu)T(Y,Y;\mu))^{1/2}$ is invariant relative to scale and orthogonal transformations, two desirable properties for measures of dependence. As a consequence this choice of measure is perhaps the most common. However, there are other choices of measures for $\mu$ that are also useful depending on the context.

\cite{dueck:2014} studied the affinely invariant distance covariance given by $\tilde T(X,Y;\mu) =T(\Sigma^{-1}_XX,\Sigma^{-1}_YY)$, where $\Sigma_X,\Sigma_Y$ are the respective covariance matrices of $X$ and $Y$ and $\mu$ is given by \eqref{eq:szikelymeas}.  They showed that the empirical version of $\tilde T(X,Y;\mu)/(\tilde T(X,Y;\mu)\tilde T(X,Y;\mu))^{1/2}$, where $\Sigma_X$ and $\Sigma_Y$ are estimated by their empirical counterparts is strongly consistent.  In addition, they provide explicit expressions in terms of special functions of the limit in the case when $X,Y$ are multivariate normal.   Further progress on this topic has been achieved in \cite{sejdinovic:sriperumbudur:gretton:fukumizu:2013} and \cite{lyons:2013}, who generalized correlation distance to a metric space.

In this paper we are interested in the empirical distance covariance and correlation applied to a stationary \seq\
$((X_t,Y_t))$ to study serial dependence, where
$X_t$ and $Y_t$ assume values in $\bbr^p$ and $\bbr^q$, respectively. We aim at an analog to the autocorrelation and autocovariance \fct s of classical \tsa\ in terms of lagged distance correlation and distance covariance.
Specifically we consider the lagged-distance covariance
\fct\
$T(X_0,Y_h;\mu)$, $h\in\bbz$, and its standardized version that takes the values in $[0,1].$  We refer to the {\em auto- and cross-distance covariances \fct s}
and their correlation analogs.  We provide \asy\ theory for the empirical auto- and cross-distance covariance and correlation
\fct s under mild conditions. Under ergodicity we prove consistency and under $\alpha$-mixing,
we derive the weak limits of the empirical  auto- and cross-distance covariance \fct s for both cases when $X_0$ and $Y_h$ are independent and dependent.  From a modeling perspective, distance correlation has limited value in providing a clear description of the nature of the dependence in the time series.  To this end, it may be difficult to find a time series model that produces a desired distance correlation.  In contrast, one could always find an autoregressive  (or more generally ARMA) process that matches the ACF for an arbitrary number of lags.  The theme in this paper will be to view the distance correlation more as a tool for testing independence rather than actually measuring dependence.

\par
The literature on distance correlation for dependent \seq s is sparse. To the best of our knowledge,
\cite{zhou:2012} was the first to study the auto-distance covariance and its empirical analog for stationary \seq s. In particular, he
proved limit theory for
$T_n(X_0,X_h;\mu)$ under so-called physical dependence \ms\ conditions
on $(X_t)$ and independence of  $X_0$ and $X_h$. \cite{fokianos:pitsillou:2015}
developed limit theory for a Ljung-Box-type statistic based on pairwise distance covariance $T_n(X_i,X_j;\mu)$
of a sample from a stationary \seq .
In both papers, the \ms\ $\mu$ is given by \eqref{eq:szikelymeas}.

 Typically, a crucial and final step in checking the quality of a fitted time series model is to examine the residuals for lack of serial dependence.  The distance correlation can be used in this regard.  However, as first pointed out in his discussion, \cite{remillard:2009} indicated that the behavior of the distance correlation when applied to the residuals of a fitted AR(1) process need not have the same limit distribution as that of the distance correlation based on the corresponding iid noise.  We provide a rigorous proof of this result for a general AR($p$) process with finite variance under certain conditions on the measure $\mu$.  Interestingly, the conditions preclude the use of the standard weight function \eqref{eq:szikelymeas} used in \cite{szekely:rizzo:bakirov:2007}.  In contrast, if the noise sequence is heavy-tailed and belongs to the domain of attraction of a stable distribution with index $\beta\in(0,2)$, the distance correlation functions for both the residuals from the fitted model and the iid noise sequence coincide.

The paper is organized as follows. In Section~\ref{sec:prelim} we commence with some basic results
for distance covariance. We give conditions on the moments of $X$ and $Y$ and the \ms\ $\mu$, which ensure that the integrals $T(X,Y;\mu)$ in \eqref{eq:i} are well-defined. We provide alternative \rep s of $T(X,Y;\mu)$
and consider various examples of finite and infinite \ms s $\mu$. Section~\ref{sec:asymp} is devoted to
the empirical auto- and cross-distance covariance and correlation \fct s. Our main results on the
\asy\ theory of these \fct s are provided in Section~\ref{subsec:asymo}. Among them are
an a.s. consistency result (Theorem~\ref{thm:consistency}) under the assumption of ergodicity and
\asy\ normality under a strong mixing condition (Theorem~\ref{thm:2}). Another main result (Theorem~\ref{thm:3}) is concerned with
the \asy\ behavior of the empirical auto-distance covariance \fct\ of the residuals of an autoregressive process for both the finite and infinite variance cases.  In Section~\ref{sec:dataex}, we
provide a small study of the empirical auto-distance correlation \fct s derived from
simulated and real-life dependent data of moderate sample size.
The proofs of  Theorems~\ref{thm:2} and
\ref{thm:3} are postponed to  Appendices~\ref{sec:proofoftheorem} and~\ref{sec:thm3}.

\section{Distance covariance for stationary time series} \label{sec:prelim}
\subsection{Conditions for existence}
From \eqref{eq:i}, the distance covariance between two vectors  $X$ and $Y$
is the squared $L^2$-distance
between the joint characteristic function of $(X,Y)$ and the product of the marginal \chf s of $X$ and $Y$
\wrt\  a \ms\ $\mu$ on $\bbr^{p+q}$.
Throughout we assume that $\mu$ is finite on sets  bounded away from the origin i.e., on sets of the form
\beam\label{eq:11}
D_\delta^c= \{(s,t): |s|\wedge |t|> \delta \}\,,\qquad \delta>0\,.
\eeam
In what follows, we interpret $(s,t)$ as a concatenated vector in $\bbr^{p+q}$ equipped with the natural norm
$|(s,t)|_{\bbr^p \times \bbr^q̈́}= \sqrt{|s|^2+ |t|^2}$. We suppress the dependence of the norm $|\cdot|$ on the dimension.
The symbol $c$ stands for any positive constant,  whose value may change from line to line, but is not of particular interest.
 Clearly if $X$ and $Y$ are independent, $T(X,Y;\mu)=0$.  On the other hand, if $\mu$ is an infinite measure, and $X$ and $Y$ are dependent,  extra conditions are needed to ensure that $T(X,Y;\mu)$ is finite.  This is the content of the following lemma.
\ble\label{lem:1}
Let $X$ and $Y$ be two possibly dependent random vectors and one of the following conditions is satisfied:
\begin{enumerate}
\item
$\mu$ is a finite \ms\ on $\bbr^{p+q}$.
\item
$\mu$ is an infinite \ms\ on $\bbr^{p+q}$, finite on the sets $D_\delta^c$, $\delta>0$, \st
\beam\label{eq:3}
\int_{\bbr^{p+q}}
 (1\wedge |s|^\alpha)\,(1\wedge |t|^\alpha) \,\mu(ds,dt) <\infty
\eeam
 and $\E[|X|^\alpha]+\E [|Y|^\alpha]<\infty$
for some $\alpha\in (0,2]$.
 \item $\mu$ is infinite  in a neighborhood of the origin and for some $\alpha\in (0,2]$, $\E[|X|^\alpha]+\E [|Y|^\alpha]<\infty$ and
\beam\label{eq:7}
\int_{\bbr^{p+q}} 1\wedge |(s,t)|^{\alpha}\,\mu(ds,dt)<\infty\,.
\eeam
\end{enumerate}
Then $T(X,Y;\mu)$ is finite.
\ele
\bre
If $\mu=\mu_1\times \mu_2$ for some \ms s $\mu_1$ and $\mu_2$ on $\bbr^p$ and $\bbr^q$, respectively, and if $\mu$ is finite on the sets
$D_\delta^c$ then it suffices for \eqref{eq:3} to verify that
\beao
\int_{|s|\le 1} |s|^\alpha\,\mu_1(ds)+ \int_{|t|\le 1} |t|^\alpha\,\mu_2(dt)<\infty\,.
\eeao
\ere
\begin{proof} (1) Since the integrand in $T(X,Y;\mu)$ is uniformly bounded the statement is trivial.\\
(2) By \eqref{eq:11},
$\mu(D_\delta^c)<\infty $ for any $\delta>0$. Therefore it remains to verify the integrability of
$|\varphi_{X,Y}(s,t)-\varphi_X(s)\,\varphi_Y(t)|^2$ on one of the sets $D_\delta$. We consider only the case $|s|\vee |t|\le 1$; the cases when
$|s|\le 1$, $|t|>1$ and $|s|>1$, $|t|\le 1$ are similar. An application of the Cauchy-Schwarz inequality yields
\beam\label{eq:6}
|\varphi_{X,Y}(s,t)-\varphi_X(s)\varphi_Y(t)|^2\le (1-|\varphi_X(s)|^2)\,(1-|\varphi_Y(t)|^2)\,.
\eeam
Since
\beao
1-|\varphi_X(s)|^2= \int_{\bbr^p} \big(1-\cos \langle s,x\rangle\big)\,\P(X-X'\in dx)
\eeao
for an independent copy $X'$ of $X$,
a Taylor expansion and the fact that $X,X'$ have finite $\alpha$th moments yield for $\alpha\in (0,2]$ and some constant $c>0$,
\beam\label{eq:10}
1-|\varphi_X(s)|^2&\le& \int_{\bbr^p}\big(
2\wedge |\langle s,x\rangle|^2\big)\,\P(X-X'\in dx)\nonumber\\
&\le &2\,\int_{|\langle s,x\rangle|\le \sqrt{2}} |\langle s,x\rangle/\sqrt{2}|^\alpha\,\P(X-X'\in dx)
+ 2\, \P(|\langle s,X-X'\rangle|>\sqrt{2}) \nonumber\\
&\le & c\,
|s|^\alpha\,\E[|X-X'|^\alpha]<\infty\,.
\eeam
In the last step we used Markov's inequality and the fact that $|\langle s,x\rangle|\le |s|\,|x|$. A corresponding bound holds for $1-|\varphi_Y(t)|^2$. Now, $T(X,Y;\mu)<\infty$ follows
from \eqref{eq:3}  and \eqref{eq:6}.\\
(3) By \eqref{eq:7}, $\mu(\{(s,t):|(s,t)|>1\})$ is finite. Therefore we need to show integrability of $|\varphi_{X,Y}(s,t)-\varphi_X(s)\varphi_Y(t)|^2$
only for $|(s,t)|\le 1$.
Using the arguments from part (2) and the finiteness of the $\alpha$th moments, we have
\beao
|\varphi_{X,Y}(s,t)-\varphi_X(s)\varphi_Y(t)|^2\le c\, (|s|^\alpha+|t|^\alpha) \le c\,|(s,t)|^{\alpha}\,.
\eeao
Now integrability of the \lhs\ at the origin  \wrt\ $\mu$ is ensured by \eqref{eq:7}.
\end{proof}
\subsection{Alternative \rep s and examples}
If $\mu=\mu_1\times \mu_2$ for \ms s $\mu_1$ and $\mu_2$ on $\bbr^p$ and $\bbr^q$ we
write for $x\in \bbr^p$ and $y\in \bbr^q$,
\beao
\hat \mu(x,y)&=& \int_{\bbr^{p+q}} \cos(\langle s,x\rangle+\langle t,y\rangle)\,\mu(ds,dt)\,,\\
\hat \mu_1(x)&=&\int_{\bbr^p} \cos \langle s,x\rangle\,\mu_1(ds)\,,\quad
\hat \mu_2(y)=\int_{\bbr^q} \cos \langle t,y\rangle\,\mu_2(dt)\,,
\eeao
for the real parts of the Fourier transforms  \wrt\ $\mu,\mu_1,\mu_2$, respectively. We assume that these transforms are well-defined.
Let $(X',Y')$ be an independent copy of $(X,Y)$, $Y'',Y'''$ are independent copies of $Y$ which are also independent of
$(X,Y), (X',Y')$. We have
\beam
T(X,Y;\mu)&=&\int_{\bbr^{p+q}} \E\Big[\ex^{i\langle s,X-X'\rangle+ i \langle t,Y-Y'\rangle} +
\ex^{i \langle s,X-X'\rangle} \,\ex^{i\langle  t, Y''-Y'''\rangle} - \ex^{i \langle s,X-X'\rangle+
i \langle t, Y-Y''\rangle}\nonumber \\&&\hspace{1.7cm}-\ex^{-i \langle s,X-X'\rangle-i \langle t ,Y-Y''\rangle}\Big]
\mu(ds,dt)\,.\label{eq:9}
\eeam
Notice that the complex-valued trigonometric \fct s under the expected value may be replaced by their real parts.
We intend to interchange the integral \wrt\ $\mu$ and the expectation.

\subsubsection{Finite $\mu$}
 For a finite \ms\ on $\bbr^{p+q}$, we may apply Fubini's theorem directly and interchange integration with expectation to obtain 
\beam
T(X,Y;\mu)
&=& \E\big[\hat \mu(X-X',Y-Y')] + \E[\hat \mu(X-X',Y''-Y''')] - 2\, \E[\hat \mu(X-X',Y-Y'')\big]\,.\nonumber\\\label{eq:13}
\eeam
If $\mu=\mu_1\times \mu_2$
we also have
\beao
T(X,Y;\mu)
&=& \E[\hat \mu_1(X-X')\,\hat \mu_2(Y-Y')] + \E[\hat \mu_1(X-X')]  \E[\hat \mu_2(Y-Y')]\\&& - 2\, \E[\hat \mu_1(X-X')\,\hat \mu_2(Y-Y'')]\,.
\eeao
\par


\subsubsection{The case of an infinite \ms\ $\mu$}
We consider an infinite  \ms\ $\mu$ on $\bbr^{p+q}$ which is finite on $D_\delta^c$ for any $\delta>0$.
We assume that $T(X,Y;\mu)$ is  finite and $\mu=\mu_1\times \mu_2$.
In this case, we cannot pass from \eqref{eq:9} to \eqref{eq:13} because the Fourier transform $\hat \mu$ is not defined as a Lebesgue integral. We have
\beam\label{eq:17}
T(X,Y;\mu)&=& 
\int_{\bbr^{p+q}} \big(\E[{\rm COS}(s,t)] +\E[{\rm
SIN}(s,t)]\big)\,\mu(ds,dt)\,,
\eeam
where
\beao
{\rm COS}(s,t)&=&\cos(\langle s,X-X'\rangle)\,\cos(\langle t,Y-Y'\rangle) +
\cos(\langle s,X-X'\rangle)\,\cos(\langle  t, Y''-Y'''\rangle)\nonumber\\&& - 2\,
\cos(\langle t,X-X'\rangle)\,\cos(\langle s, Y-Y''\rangle)\,,\\
{\rm SIN}(s,t)&=&-\sin(\langle s,X-X'\rangle)\,\sin(\langle t,Y-Y'\rangle) -
\sin(\langle s,X-X'\rangle)\,\sin(\langle  t, Y''-Y'''\rangle)\nonumber\\&& + 2\,
\sin(\langle t,X-X'\rangle)\,\sin(\langle s, Y-Y''\rangle)\,.
\eeao
Using the fact that
\beao
\cos u\,\cos v= 1-(1-\cos u)-(1-\cos v)+ (1-\cos u)(1-\cos v)\,,
\eeao
calculation shows that
\beao
\E\big[{\rm COS}(s,t)]&=&
\E\big[
(1-\cos(\langle s,X-X'\rangle))\,(1-\cos(\langle t,Y-Y'\rangle))\nonumber\\&& +
(1-\cos(\langle s,X-X'\rangle))\,(1-\cos(\langle  t, Y''-Y'''\rangle))\nonumber\\&& - 2\,(1-
\cos(\langle t,X-X'\rangle))\,(1-\cos(\langle s, Y-Y''\rangle))\big]\,. 
\eeao
A Taylor series argument shows that for $\alpha\in(0 ,2]$,
\beao
\E[|{\rm COS}(s,t)|]
&\le& c\,\Big(
\E\big[(1\wedge |\langle s,X-X'\rangle/\sqrt{2}|^\alpha)\,(1\wedge |\langle t,Y-Y'\rangle/\sqrt{2}|^\alpha)\big]\\
&&+\E\big[1\wedge |\langle s,X-X'\rangle/\sqrt{2}|^\alpha\big]\,\E\big[1\wedge |\langle t,Y-Y'\rangle/\sqrt{2}|^\alpha\big]\\
&&+\E\big[(1 \wedge |\langle t,X-X'\rangle/\sqrt{2}|^\alpha)\,(1\wedge |\langle s, Y-Y''\rangle/\sqrt{2}|^\alpha\big|\big]\Big)\,.
\eeao
Under condition \eqref{eq:3}
the \rhs\ is integrable \wrt\ $\mu$ if
\beam\label{eq;3a}
\E[|X|^\alpha+|Y|^\alpha+ |X|^\alpha\,|Y|^\alpha]<\infty\,.
\eeam
An application of Fubini's theorem yields
\beao
\int_{\bbr^{p+q}} \E[{\rm COS}(s,t)]\,\mu(ds,dt)
&=& \E\Big[\int_{\bbr^{p+q}}\Big( (1-\cos(\langle s,X-X'\rangle))\,(1-\cos(\langle t,Y-Y'\rangle))\nonumber\\&& +
(1-\cos(\langle s,X-X'\rangle))\,(1-\cos(\langle  t, Y''-Y'''\rangle))\nonumber\\&& - 2\,(1-
\cos(\langle t,X-X'\rangle))\,(1-\cos(\langle s, Y-Y''\rangle))\Big)\,\mu(ds,dt)\Big]\,.
\eeao
If we assume that the restrictions $\mu_1,\mu_2$ of $\mu$ to $\bbr^p$ and $\bbr^q$ are symmetric about the origin then we have
\beao
 \E\big[{\rm SIN}(s,t)\big]= -\E\big[{\rm SIN}(-s,t)\big]=-\E\big[{\rm SIN}(s,-t)\big]\,.
\eeao
Together with the symmetry property of $\mu$ this implies that $\int_{\bbr^{p+q}}\E\big[{\rm SIN}(s,t)\big]\,\mu(ds,dt)=0$.
\par
We summarize these arguments. For any \ms\ $\nu$ on $\bbr^d$ we write
\beao
\tilde \nu(s)= \int_{\bbr^d}(1-\cos \langle s,x\rangle) \,\nu(dx)\,,\qquad s\in\bbr^d\,.
\eeao
\ble Assume \eqref{eq:3} and \eqref{eq;3a} for some $\alpha\in (0,2]$. If $\mu_1,\mu_2$ are symmetric about the origin and $\mu=\mu_1\times \mu_2$
then
\beam\label{eq:14}
T(X,Y;\mu) &=& \E[\tilde \mu_1 (X-X')\,\tilde\mu_2(Y-Y')] +\E[\tilde \mu_1(X-X')]\,\E[ \tilde \mu_2(Y-Y')]\nonumber\\&&- 2\,\E[\tilde\mu_1(X-X')
\tilde \mu_2(Y-Y'')]\,.
\eeam
\ele
\bre
For further use, we mention the alternative \rep\ of \eqref{eq:14}:
\beam\label{eq:covex}
T(X,Y;\mu)&=&\cov\big(\tilde \mu_1 (X-X'),\tilde\mu_2(Y-Y')\big) -2\,\cov\big(\E[\tilde \mu_1 (X-X')\mid X]\,,\E[\tilde\mu_2(Y-Y')\mid Y]
\big)\,.\nonumber\\
\eeam
\ere
\bexam\label{exam:szekly}
Assume that $\mu$ has density $w$ on $\bbr^{p+q}$ given by
\beam\label{eq:szeklyi}
w(s,t)= c_{p,q}\,|s|^{-\alpha-p}\,|t|^{-\alpha-q}\,,\qquad s\in\bbr^p\,,t\in\bbr^q\,,
\eeam
for some positive constant $c_{p,q} = c_pc_q$. For any $d\ge 1$ and $\alpha\in (0,2)$,  one can choose $c_{d}$ \st\
\beam\label{eq:density}
\int_{\bbr^{d}} (1- \cos\langle s,x\rangle)\,c_d\,|s|^{-\alpha-d}\,ds= |x|^{\alpha}\,.
\eeam
Under the additional moment assumption \eqref{eq;3a}
we obtain from  \eqref{eq:14}
\beam\label{eq:15}
T(X,Y;\mu)&=&\E[|X-X'|^\alpha\,|Y-Y'|^\alpha] + \E[|X-X'|^\alpha] \, \E[Y-Y'|^\alpha] - 2\, \E[|X-X'|^\alpha\,|Y-Y''|^\alpha]\,.\nonumber\\
\eeam
This is the {\em distance covariance} introduced by
\cite{szekely:rizzo:bakirov:2007}.
\eexam
The distance covariance $T(X,Y;\mu)$ introduced in \eqref{eq:15} has several good properties.
 It is homogeneous under positive scaling and is also invariant under orthonormal
transformations of $X$ and $Y$. Some of these properties are shared with other distance covariances
when $\mu$ is infinite. We illustrate this for a \levy\ \ms\ $\mu$ on $\bbr^{p+q}$, i.e., it satisfies \eqref{eq:7} for $\alpha=2$.
In particular, $\mu$ is finite on sets bounded away from zero. Via the \levy -Khintchine formula,
a \levy\ \ms\ $\mu$  corresponds
 to an $\bbr^{p+q}$-valued
\id\ random vector $(Z_1,Z_2)$ (with $Z_1$ assuming values in $\bbr^p$ and $Z_2$ in $\bbr^q$) and \chf
\beam\label{eq:00}
\varphi_{Z_1,Z_2}(x,y) =\exp\Big\{-\int_{\bbr^{p+q}}
\Big(\ex^{i\langle s , x\rangle+i\langle t,y\rangle}- 1-(i\langle x,s\rangle+i\langle y,t\rangle)\1 (|(s,t)|\le 1) \Big)\,\mu(ds,dt)\Big\}\,.\nonumber\\
\eeam
\ble
Assume that there exists an $\alpha\in (0,2]$ \st\ $\E[|X|^\alpha]+\E[|Y|^\alpha]<\infty$ and $\mu$ is a symmetric \levy\ \ms\ corresponding to \eqref{eq:00}
\st\ \eqref{eq:7} holds.
Then
\beam\label{eq:18}
T(X,Y;\mu)&=& {\rm Re}\,\E \Big[-\log \varphi_{Z_1,Z_2} (X-X',Y-Y')-
\log \varphi_{Z_1,Z_2}(X-X',Y''-Y''')\nonumber\\&&+2
\log \varphi_{Z_1,Z_2}(X-X', Y-Y'')\Big]\,.
\eeam
\ele
\bre
We observe that \eqref{eq:18} always vanishes if $Z_1$ and $Z_2$ are independent.
\ere

\begin{proof}
By the symmetry of the random vectors in \eqref{eq:9} and the measure $\mu$, we have
\beao
\lefteqn{{\rm Re}\int_{\bbr^{p+q}}\E\Big[\ex^{i\langle s,X-X'\rangle+ i \langle t,Y-Y'\rangle}-1\Big] \,\mu(ds,dt)}\\&=&
{\rm Re} \int_{\bbr^{p+q}}\E\Big[\ex^{i\langle s,X-X'\rangle+ i \langle t,Y-Y'\rangle} -
1-(i\langle s,X-X'\rangle+ i \langle t,Y-Y'\rangle)\,\1\big(|(s,t)|\le 1\big)\Big]\,\mu(ds,dt)\\
&=&  {\rm Re}\,\E\big[-\log \varphi_{Z_1,Z_2}(X-X',Y-Y')\big]\,.
\eeao
The last step is justified if we can interchange the integral and the expected value.
Therefore we have to verify that the following integral is finite:
\beao
\int_{\bbr^{p+q}}\E\Big[\Big|\ex^{i\langle s,X-X'\rangle+ i \langle t,Y-Y'\rangle} -
1-(i\langle s,X-X'\rangle+ i \langle t,Y-Y'\rangle)\,\1\big(|(s,t)|\le 1\big)\Big|\Big]\,\mu(ds,dt)
\,.
\eeao
We denote the integrals over the disjoint sets $\{(s,t):|(s,t)|\le 1\}$ and $\{(s,t): |(s,t)|>1\}$ by $I_1$ and $I_2$, respectively.
The quantity $I_2$ is bounded since the integrand is bounded and $\mu$ is finite on sets bounded away from zero. A Taylor expansion
shows for $\alpha\in (0,2]$,
\beao
I_1&\le& c\,\int_{|(s,t)|\le 1} \E\big[2\wedge  (|\langle s ,X-X'\rangle|+|\langle t,Y-Y'\rangle|)^2\big] \,\mu(ds,dt)\\
&\le& c\, (\E|X|^\alpha]+\E|Y|^\alpha])\,\int_{|(s,t)|\le 1} 1\wedge  |(s,t)|^\alpha  \,\mu(ds,dt)\\
\eeao
and the \rhs\ is finite by assumption.
\par
Proceeding in the same way as above for the remaining expressions in \eqref{eq:9}, the lemma is proved.
\end{proof}

\bexam
Assume that $\mu$ is a \pro y \ms\ of a random vector $(Z_1,Z_2)$ in $\bbr^{p+q}$ and that $Z_1$ and $Z_2$ are independent.
Then
\beao
T(X,Y;\mu)&=&\E [\varphi_{Z_1}(X-X')\,\varphi_{Z_2}(Y-Y')] + \E[\varphi_{Z_1}(X-X')]\,\E[\varphi_{Z_2}(Y''-Y''')]\\&& - 2\, \E[\varphi_{Z_1}(X-X')\,\varphi_{Z_2}(Y-Y'')]\,.
\eeao
For example, consider independent symmetric $Z_1$ and $Z_2$ with multivariate $\beta$-stable \ds s in $\bbr^p$ and $\bbr^q$,
respectively, for some $\beta\in (0,2]$.  
They have joint  \chf\ given by
$\varphi_{Z_1,Z_2}(x,y)= \ex^{-(|x|^\beta+|y|^\beta)}$. Therefore
\beam\label{eq:betastab}
T(X,Y;\mu)&=&\E [\ex^{-(|X-X'|^\beta+|Y-Y'|^\beta)}] + \E[\ex^{-|X-X'|^\beta}]\,\E[\ex^{-|Y-Y'|^\beta}]- 2\,
\E [\ex^{-(|X-X'|^\beta+|Y-Y''|^\beta)}]\,.\nonumber\\
\eeam
\eexam

\bexam
Assume that $X$ and $Y$ are integer-valued. Consider the spectral densities $w_1$ and $w_2$ on $[-\pi,\pi]$ of
two real-valued second-order stationary processes and assume $\mu(s,t)=w_1(s)w_2(t)$. Denote the covariance \fct s on the integers
corresponding to $w_1$ and $w_2$ by $\gamma_1$ and $\gamma_2$, respectively. We have the well-known relation
\beao
\int_{-\pi}^\pi \ex^{itk} \,w_i(t)\,dt=\int_{-\pi}^\pi \cos(tk)\, \,w_i(t)\,dt
=\gamma_i(k)\,,\qquad k\in\bbz\,,
 \eeao
where we also exploit the symmetry of the \fct s $w_i$.
If we restrict integration in \eqref{eq:9} to $[-\pi,\pi]^2$ we obtain, abusing notation,
\beao
T(X,Y;\mu)&=&
\E[\gamma_1(X-X')\,\gamma_2(Y-Y')] + \E[\gamma_1(X-X')]  \E[\gamma_2(Y-Y')]\\&& - 2\, \E[\gamma_1(X-X')\,\gamma_2(Y-Y'')]\,.
\eeao
The spectral density of a stationary process may have singularities (e.g.\ for fractional ARMA processes)
but this density is integrable on $[-\pi,\pi]$. If $w_1,w_2$ are positive Lebesgue a.e.\ on $[0,\pi]$ then $T(X,Y;\mu)=0$ \fif\ $X,Y$ are
independent. Indeed, the \chf\ of an integer-valued \rv\ is periodic with period $2\pi$.
\eexam

\bexam
To illustrate  \eqref{eq:18} we consider a symmetric $\alpha$-stable vector $(Z_1,Z_2)$ for $\alpha\in (0,2)$ with log-\chf\
\beao
-\log \varphi_{Z_1,Z_2}(x,y)= \int _{\bbs^{p+q-1}} |\langle s,x\rangle+\langle t,y\rangle |^\alpha \,m(ds,dt)
\eeao
and $m$ is a finite symmetric \ms\ on the unit sphere $\bbs^{p+q-1}$ of $\bbr^{p+q}$. Then we have
\beao
T(X,Y;\mu) &=& \int _{\bbs^{p+q-1}} \E\big[|\langle s,X-X'\rangle+\langle t,Y-Y'\rangle|^\alpha+ |\langle s,X-X'\rangle+\langle t,Y''-Y'''\rangle |^\alpha
\\&&-2 \,|\langle s,X-X'\rangle +\langle t,Y'-Y''\rangle|^\alpha\big] \,m(ds,dt)\,.
\eeao
A special case is the sub-Gaussian $\alpha/2$-stable random vectors with \chf\
\beao
-\log \varphi_{Z_1,Z_2}(x,y)=  |(x,y)'\Sigma (x,y)|^{\alpha/2}\,,
\eeao
where $\Sigma$ is the covariance matrix of an $\bbr^{p+q}$-valued random vector and we write $(x,y)$ for the concatanation of any $x\in \bbr^p$ and
$y\in\bbr^q$.
Then
\beao
T(X,Y;\mu) &=& \E\big[|(X-X',Y-Y')'\Sigma\,(X-X',Y-Y')|^{\alpha/2}\\&&+[|(X-X',Y''-Y''')'\Sigma\,(X-X',Y''-Y''')|^{\alpha/2}
\\&&-2 [|(X-X',Y-Y'')'\Sigma\,(X-X',Y-Y'')|^{\alpha/2}\big]\,.
\eeao
In particular, if $\Sigma$ is block-diagonal with $\Sigma_1$ a $p\times p$ covariance matrix and
$\Sigma_2$ a $q\times q$ covariance matrix, we have
\beao
T(X,Y;\mu) &=& \E\big[|(X-X')'\Sigma_1\,(X-X') + (Y-Y')'\Sigma_2(Y-Y')|^{\alpha/2}\\&&
+|(X-X')'\Sigma_1\,(X-X') + (Y''-Y''')'\Sigma_2(Y''-Y''')|^{\alpha/2}\\&&
-2 |(X-X')'\Sigma_1\,(X-X') + (Y-Y'')'\Sigma_2(Y-Y'')|^{\alpha/2}\big]\,,
\eeao
and if $\Sigma$ is the identity matrix,
\beam\label{eq:20}
T(X,Y;\mu) &=& \E\big[\big|\,|X-X'|^2+|Y-Y'|^2\big|^{\alpha/2}+\big|\,|X-X'|^2+|Y''-Y'''|^2\big|^{\alpha/2}\nonumber\\
&&-2\big|\,|X-X'|^2+|Y-Y''|^2\big|^{\alpha/2}\big]\,.
\eeam
We notice that for these examples, $T(X,Y;\mu)$ is scale homogeneous 
($T(cX,cY;\mu)=|c|^\alpha T(X,Y;\mu)$) and \eqref{eq:20} is invariant 
under orthonormal transformations ($T(RX,SY;\mu)=T(X,Y;\mu)$ for orthonormal matrices $R$ and $S$), 
properties also enjoyed by the weight function in Example~\ref{exam:szekly}.
\eexam

\section{The empirical distance covariance \fct\ of a stationary \seq }\label{sec:asymp}
In this section we consider the empirical distance covariance for a stationary time series $((X_t,Y_t))$ with
generic element $(X,Y)$ where $X$ and $Y$ assume values in $\bbr^p$ and $\bbr^q$, respectively. The empirical distance
covariance is given by
\beao
T_n(X,Y;\mu) &=& \int_{\bbr^{p+q}} \big|\varphi_{X,Y}^n(s,t)-\varphi_X^n(s)\,\varphi_Y^n(t)\big|^2\,\mu(ds,dt)\,,
\eeao
where the empirical \chf\ is given by
\beao
\varphi_{X,Y}^n(s,t)= \dfrac 1 n \sum_{j=1}^n \ex^{i\,\langle s,X_j\rangle+i\,\langle t,Y_j\rangle}\,,\qquad n\ge 1\,,
\eeao
and $\varphi_X^n(s)= \varphi_{X,Y}^n(s,0)$ and $\varphi_Y^n(s)= \varphi_{X,Y}^n(0,t)$.
\subsection{Asymptotic results for the empirical distance correlation}\label{subsec:asymo}
 Under the conditions of  Lemma \ref{lem:1} that ensure the finiteness of $T(X,Y;\mu)$, we show that  $T_n$ is consistent for stationary ergodic time series.
\bth \label{thm:consistency}
 Consider a stationary ergodic time series
$((X_i,Y_i))_{i=1,2,\ldots}$ with values in $\bbr ^{p+q}$ and assume one of the three conditions in Lemma \ref{lem:1} are satisfied.  Then
$$
T_n(X,Y;\mu)\stas T(X,Y;\mu)\,~~~\mbox{as}~~~\nto.
$$
\ethe
\begin{proof}
We denote the differences of the \chf s and its empirical analog  by
\beao
C(s,t)&=& \varphi_{X,Y}(s,t)-\varphi_X(s)\varphi_Y(t)\,,\\
C_n(s,t)&=&\varphi^n_{X,Y}(s,t)-\varphi^n_X(s)\varphi^n_Y(t)\,,\qquad (s,t)\in \bbr^{p+q}\,.
\eeao
Each of the processes $\varphi^n_{X,Y}$, $\varphi^n_X$, $\varphi^n_Y$ is a sample mean of iid bounded continuous processes defined
on $\bbr^{p+q}$. Consider the compact set
\beam\label{eq:3w}
K_\delta=\{(s,t)\in\bbr^{p+q}: \delta\le  |s|\wedge |t|\,, |s|\vee |t|\le 1/\delta\}
\eeam
for small $\delta>0$. By the ergodic theorem on
${\mathcal C}( K_\delta)$, the space of continuous \fct s on $K_\delta$,
$\varphi^n_{X,Y}\stas \varphi_{X,Y}$ as $\nto$; see \cite{krengel:1985}. Hence
\beao
\int_{K_\delta} |C_n(s,t)|^2\,\mu(ds,dt)\stas \int_{K_\delta} |C(s,t)|^2\,\mu(ds,dt)\,,\qquad\nto \,.
\eeao
It remains to show that
\beao
\lim_{\delta \downarrow 0}\limsup_{\nto} \int_{K_\delta^c} |C_n(s,t)|^2\,\mu(ds,dt)=0\quad \as
\eeao
If $\mu$ is a finite \ms\ we have
\beao
\lim_{\delta\downarrow 0} \limsup_{\nto}\int_{K_\delta^c}|C_n(s,t)|^2 \mu(ds,dt)\le c\,\lim_{\delta \downarrow 0} \mu(K_\delta^c)=0\,.
\eeao
Now assume that $\mu$ is infinite on the axes or at zero and \eqref{eq:3} holds.
We apply inequality \eqref{eq:6} under the assumption that $(X,Y)$ has the empirical \pro y \ms\ of the sample $(X_i,Y_i)$, $i=1,\ldots,n$.
Since the empirical \ms\ has all
moments finite we obtain from \eqref{eq:10} that for $\alpha\in (0,2]$,
\beao
1- |\varphi_{X}^n(s)|^2\le c\,|s|^\alpha\,\E_{n,X} [|X-X'|^\alpha]= c\, |s|^\alpha\,n^{-2} \sum_{1\le k,l\le n}|X_k-X_l|^\alpha\,,
\eeao
where $X,X'$ are independent and each of them has the empirical \ds\ of the $X$-sample. The \rhs\ is a $U$-statistic which converges
a.s. to $\E[|X-X'|^\alpha]$ as $\nto$ provided this moment is finite. This follows from the ergodic theorem for $U$-statistics; see
\cite{aaronson:dehling}.
The same argument as for
part (2) of Lemma~\ref{lem:1} implies that on $K_\delta^c$,
\beao
|C_n(s,t)|^2\le c\, \E_{n,X} [|X-X'|^\alpha]\,\E_{n,Y} [|Y-Y'|^\alpha]\, (1\wedge |s|^\alpha)\,(1\wedge |t|^\alpha)\,.
\eeao
By the ergodic theorem,
\beao
\limsup_{\nto} \int_{K_\delta^c}|C_n(s,t)|^2 \mu(ds,dt)\le c\,\E[|X-X'|^\alpha]\,\E[|Y-Y'|^\alpha]\,\int_{K_\delta^c}
(1\wedge |s|^\alpha)(1\wedge |t|^\alpha) \mu(ds,dt)\quad \as
\eeao
and the latter integral converges to zero as $\delta\downarrow 0$ by assumption.
\par
If the \ms\ $\mu$ is infinite at zero and \eqref{eq:7} holds the proof is analogous.
\end{proof}

In order to prove weak \con\ of $T_n$ we assume that the \seq\ $((X_i,Y_i))$ with values in $\bbr^{p+q}$ is
$\alpha$-mixing with rate \fct\ $(\alpha_h)$; see \cite{ibragimov:linnik:1971} and \cite{doukhan:1994}
for the corresponding definitions.
 We have the following result.
\bth \label{thm:2}
Assume that $((X_i,Y_i))$ is a strictly stationary \seq\ with values in $\bbr^{p+q}$ \st\
$\sum_h \alpha_h^{1/r}<\infty$ for some $r>1$. Set $u=2r/(r-1)$ and write $X=(X^{(1)},\ldots, X^{(p)})$ and $Y=(Y^{(1)},\ldots, Y^{(q)})$.
\begin{enumerate}
\item
Assume that $X_0$ and $Y_0$ are independent and for some $\alpha\in (u/2,u]$, $\epsilon\in [0,1/2)$ and $\alpha'\le \min(2,\alpha)$, the following hold:
\beam\label{eq:moment}
\E[|X|^\alpha+|Y|^\alpha]<\infty,\qquad~~\qquad \E \big[\prod_{l=1}
^p |X^{(l)}|^\alpha\big]<\infty\,,\quad\E \big[\prod_{l=1}^q |Y^{(l)}|^\alpha\big]<\infty\,,
\eeam
and
\beam\label{eq:j}
\int_{\bbr^{p+q}} (1\wedge |s|^{\alpha'(1+\epsilon)/u})(1\wedge |t|^{\alpha'(1+\epsilon)/u})\,\mu(ds,dt)<\infty\,.
\eeam
Then
\beam\label{eq:t_n}
n\,T_n(X,Y;\mu)\std \|G\|_\mu^2=\int_{\bbr^{p+q}} |G(s,t)|^2\,\mu(ds,dt)\,,
\eeam
where $G$ is a complex-valued mean-zero Gaussian process whose covariance structure is given in \eqref{eq:covG} with $h=0$ and depends on the dependence structure of $((X_t,Y_t))$.
\item
Assume that $X_0$ and $Y_0$ are dependent and for some $\alpha\in
     (u/2,u]$, $\epsilon\in [0,1/2)$ and for $\alpha'\le \min(2,\alpha)$ the following hold:
\beam
\label{comd:mom}
\E[|X|^{2\alpha}+|Y|^{2\alpha}]<\infty,\quad \E\Big[ (1\vee \prod_{l=1}^p |X^{(l)}|^\alpha)(1\vee \prod_{k=1}^q
|Y^{(k)}|^\alpha)\Big]<\infty\,,
\eeam
and
\beam\label{eq:j1}
\int_{\bbr^{p+q}} (1\wedge |s|^{\alpha'(1+\epsilon)/u})(1\wedge |t|^{\alpha'(1+\epsilon)/u})\,\mu(ds,dt)<\infty\,.
\eeam
Then
\beam\label{eq:t_n1}
\sqrt{n}\,(T_n(X,Y;\mu)-T(X,Y;\mu)) \std G_\mu' =\int_{\bbr^{p+q}} G'(s,t)\,\mu(ds,dt)\,,
\eeam
where $G'(s,t)=2{\rm Re}\{G(s,t)C(s,t)\}$ is a mean-zero Gaussian process.
\end{enumerate}
\ethe
The proof of Theorem~\ref{thm:2} is given in Appendix~\ref{sec:proofoftheorem}.
\bre
 We notice that \eqref{eq:j} and \eqref{eq:j1} are always satisfied if $\mu$ is a finite \ms.
\ere
\bre\label{rem:iid}
If $(X_i)$ and $(Y_i)$ are two independent iid \seq s then the statement of Theorem~\ref{thm:2}(1)
remains valid if for some $\alpha\in (0,2]$,  $\E [|X|^\alpha]+ \E [|Y|^\alpha]<\infty$ and
\beam\label{eq:1958}
\int_{\bbr^{p+q}} (1\wedge |s|^{\alpha})(1\wedge |t|^{\alpha})\,\mu(ds,dt)<\infty\,.
\eeam
\ere
\bre
The \ds\ of the limit variable in \eqref{eq:t_n} is generally not tractable. Therefore one must use numerical or resampling methods for determining quantiles of $nT_n(X,Y;\mu)$.  On the other hand, the limit distribution in  \eqref{eq:t_n1} is normally distributed with mean 0 and covariance $\sigma^2_\mu$ that can be easily calculated from the covariance function of $G(s,t)$ and $C(s,t)$.  Notice that if $C(s,t)=0$, the limit random variable in \eqref{eq:t_n1} is 0 and part (1) of the theorem applies.  Again resampling or subsampling methods must be employed to determine quantiles of $nT_n$.
\ere

\subsection{Testing serial dependence for multivariate time series}
Define the {\em cross-distance covariance \fct } (CDCVF) of a strictly stationary \seq\ $((X_t,Y_t))$ by
\beao
T^{X,Y}_\mu(h)=T(X_0,Y_h;\mu)\,, \qquad h\in\bbz\,,
\eeao
and the {\em auto-distance covariance \fct } (ADCVF) of a stationary \seq\ $(X_t)$
by
\beao
T_\mu^X(h)=T_\mu^{X,X}(h)\,, \qquad h\in\bbz\,.
\eeao
Here and in what follows, we assume that $\mu=\mu_1\times \mu_2$ for suitable \ms s $\mu_1$ on $\bbr^p$ and $\mu_2$ on $\bbr^q$. In the case of an ADCVF we also assume $\mu_1=\mu_2$.
The empirical versions $T_{n,\mu}^X$ and $T_{n,\mu}^{X,Y}$ are defined correspondingly.
For example, for integer $h\ge 0$, one needs to replace $\varphi_{X,Y}^n(s,t)$ in the definition of $T_n(X,Y;\mu)$ by
\beao
\varphi_{X_0,Y_h}^n(s,t) = \frac 1 {n} \sum_{j=1}^{n-h}\ex^{i \,\langle s,X_j\rangle +i\, \langle t,Y_{j+h}\rangle},\qquad s\in\bbr^p\,,t\in\bbr^q\,, \quad  n\ge h+1\,,
\eeao
with the corresponding modifications for the marginal empirical \chf s. For finite $h$, the change from the upper summation limit
$n$ to $n-h$ has no influence on the \asy\ theory.
\par
We also introduce the corresponding {\em cross-distance correlation \fct } (CDCF)
\beao
R^{X,Y}_\mu(h)=\dfrac{T_\mu^{X,Y}(h)}{\sqrt{T_\mu^X(0)\, T_\mu^Y(0)}}\,, \qquad h\in\bbz\,,
\eeao
and
 {\em auto-distance correlation \fct } (ADCF)
\beao
R^{X}_\mu(h)=\dfrac{T_\mu^{X}(h)}{T_\mu^X(0)}\,, \qquad h\in\bbz\,,
\eeao
The quantities $R^{X,Y}_\mu(h)$ assume values in $[0,1]$, with the two endpoints representing independence and complete dependence.
The empirical CDCF $R_{n,\mu}^{X.Y}$ and ADCF $R_{n,\mu}^X$ are defined by replacing the distance covariances $T_\mu^{X,Y}(h)$ by the corresponding empirical versions
$T_{n,\mu}^{X,Y}(h)$.
\par
The empirical ADCV was examined in \cite{zhou:2012} and \cite{fokianos:pitsillou:2015} as an alternative
tool for testing serial dependence, in the way that it also captures non-linear dependence.
They always choose the \ms\ $\mu=\mu_1\times \mu_1$ with density \eqref{eq:szeklyi}.
\par
In contrast to the autocorrelation and cross-correlation \fct s of standard stationary \ts\ models (such as ARMA, GARCH)
it is in general complicated (or impossible)
to provide explicit (and tractable) expressions for $T_\mu^X(h)$ and $T^{X,Y}_\mu(h)$ or even to say anything about the
rate of decay of these quantities when $h\to \infty$. However, in view of \eqref{eq:covex} we observe that
\beao
T_\mu^X(h) = \cov\big(\tilde \mu_1 (X_0-X_0'),\tilde\mu_1(X_h-X_h')\big) -
2\,\cov\big(\E[\tilde \mu_1 (X_0-X_0')\mid X_0]\,,\E[\tilde\mu_1(X_h-X'_0)\mid X_h]
\big)\,.
\eeao
 While this is not the autocovariance \fct\ of a stationary process, it is possible to bound each of the terms in case  $(X_t)$ is $\alpha$-mixing with rate \fct\ $(\alpha_h)$.  In this case, one may use bounds for the autocovariance \fct s of the
stationary series $(\tilde\mu_1(X_t-X_t'))$ and $(\E[\tilde\mu_1(X_t-X'_0)\mid X_t])$ which inherit $\alpha$-mixing from $(X_t)$ with the same rate \fct .
For example, a
standard inequality (\cite{doukhan:1994}, Section 1.2.2, Theorem 3(a)) yields that
\beao
T^X_\mu(h)\le c\,\alpha_h^{1/r}\, \big(\E[(\tilde \mu_1 (X_0-X_0'))^{u}]\big)^{2/u}
\eeao
for positive $c$ and $r>0$ \st\ $r^{-1}+ 2u^{-1}=1$. If $\tilde \mu_1$ is bounded we also have $T^X_\mu(h)\le c\,\alpha_h$ for some positive constant.
Similar bounds can be found for $T^{X,Y}_\mu(h)$ provided $((X_t,Y_t))$ is $\alpha$-mixing.
\par
Next we give an example where the ADCVF can be calculated explicitly.
\bexam
Consider a univariate strictly stationary Gaussian \ts\ $(X_t)$ with
mean zero, variance $\sigma^2$  and autocovariance \fct\ $\gamma_X$.
We choose
a  Gaussian \pro y \ms\ $\mu$ which leads to the relation \eqref{eq:betastab}. Choose $N_1,N_2,N_3$ iid $N(0,2)$-distributed
independent of the independent quantities $(X_0,X_h), (X_0',X_h'), X_h''$. Then for $h\ge 0$,
\beao
T_\mu^X(h)&=&\E\big[\ex^{i\,N_1(X_0-X_0')+ iN_2 (X_h-X_h')}\big]+ \big(\E\big[\ex^{i\,N_1(X_0-X_0')}\big]\big)^2
-2 \,\E\big[\ex^{i\,N_1(X_0-X_0')+ iN_2 (X_h-X_h'')}\big] \\
&=& \E\big[\ex^{i\,(N_1X_0+N_2 X_h)-i (N_1X_0'+N_2X_h')}\big]+\big(\E\big[\ex^{i\,N_1(X_0-X_0')}\big]\big)^2\\
&&-2\,\E\big[\ex^{i\,(N_1X_0+N_2 X_h)-i(N_1X_0'+N_2X_h'')}\big] \\
&=&\E\big[\ex^{i\,N_3 \big(N_1^2\sigma^2+N_2^2 \sigma^2 +2 \gamma_X(h)N_1N_2\big)^{1/2}}\big]
+ \big(\E\big[\ex^{i\,N_3 (N_1^2\sigma^2)^{1/2}}\big]\big)^2\\
&&-2\,\E\big[\ex^{i\,N_3\, \big(N_1^2\sigma^2+N_2^2 \sigma^2 +\gamma_X(h)N_1N_2\big)^{1/2}}\big]\\
&=&\E\big[\ex^{-\big(N_1^2\sigma^2+N_2^2 \sigma^2 +2 \gamma_X(h)N_1N_2\big)}\big]
+ \big(\E\big[\ex^{- N_1^2\sigma^2}\big]\big)^2\\
&&-2\,\E\big[\ex^{-\big(N_1^2\sigma^2+N_2^2 \sigma^2 +\gamma_X(h)N_1N_2\big)}\big]\,.
\eeao
For the evaluation of this expression we focus on the first term, the other cases being similar.
Observing that $\sigma^2\pm \gamma_X(h)$ are the eigenvalues of the covariance matrix
\beao
\left( \barr{cc}\sigma^2& \gamma_X(h)\\
\gamma_X(h)&\sigma^2 \earr\right)\,,
\eeao
calculation shows that
\beao
N_1^2\sigma^2+N_2^2 \sigma^2 +2 \gamma_X(h)N_1N_2\eqd N_1^2  (\sigma^2-\gamma_X(h)) + N_2^2 (\sigma^2+\gamma_X(h))\,.
\eeao
Now the moment generating \fct\ of a $\chi^2$-distributed \rv\  yields
\beao
\E\big[\ex^{-\big(N_1^2\sigma^2+N_2^2 \sigma^2 +2 \gamma_X(h)N_1N_2\big)}\big]
= \big(1+ 4 (\sigma^2-\gamma_X(h))\big)^{-1/2}\,\big(1+ 4 (\sigma^2+\gamma_X(h))\big)^{-1/2}\,.
\eeao
Proceeding in a similar fashion, we obtain
\beao
T_\mu^X(h)&=&\big(1+ 4 (\sigma^2-\gamma_X(h))\big)^{-1/2}\,\big(1+ 4 (\sigma^2+\gamma_X(h))\big)^{-1/2}
+(1+ 4 \sigma^2)^{-1}\\
&&- 2\, \big(1+ 4 (\sigma^2-\gamma_X(h)/2)\big)^{-1/2}\,\big(1+ 4 (\sigma^2+\gamma_X(h)/2)\big)^{-1/2}\,.
\eeao
If $\gamma_X(h)\to 0$ as $h\to\infty$ Taylor expansions yield
$T_\mu^X(h) \sim 4\gamma_X^2(h)/(1+ 4\sigma^2)^3$.
A similar result was given in \cite{fokianos:pitsillou:2015}, where they derived an explicit expression for $T_\mu^X(h)$ for a stationary Gaussian process $(X_t)$ with weight function \eqref{eq:szikelymeas}.
\eexam
If $((X_t,Y_t))$ is strictly stationary and ergodic then
$((X_t,Y_{t+h}))$ is a strictly stationary ergodic \seq\ for every integer $h$. Then Theorem~\ref{thm:consistency} applies.
\bco Under the conditions of Theorem~\ref{thm:consistency}, for $h\ge 0$,
\beao
T_{n,\mu}^{X,Y}(h)\stas T_\mu^{X,Y}(h)\qquad\mbox{and}\qquad  T_{n,\mu}^X(h)\stas T_\mu^X(h)\,,
\eeao
and
\beao
R_{n,\mu}^{X,Y}(h)\stas R_\mu^{X,Y}(h)\qquad\mbox{and}\qquad  R_{n,\mu}^X(h)\stas R_\mu^X(h)\,.
\eeao
\eco
Applying Theorem \ref{thm:2} and Theorem~\ref{thm:consistency}, we also have the
following weak dependence result under $\alpha$-mixing.
\cite{zhou:2012} proved the corresponding
result under conditions on the so-called {\em physical dependence \ms }.
\bco\label{cor:adcv}
Assume that $X_0$ and $Y_h$ are independent for some $h\ge 0$ and the \seq\ $((X_t,Y_t))$ satisfies the conditions of Theorem~\ref{thm:2}.
Then
\beao
n\,T_{n,\mu}^{X,Y}(h)\std \|G_h\|_\mu^2\qquad\mbox{and}\qquad n\,R_{n,\mu}^{X,Y}(h)\std \dfrac{\|G_h\|_\mu^2}{
\sqrt{T_\mu^X(0)\,T_\mu^Y(0)}}\,,
\eeao
where $G_h$ is a centered Gaussian process on $\bbr^{p+q}$.
\eco
\bre\label{rem:5k} From the proof of Theorem~\ref{thm:2} (the central limit theorem for the multivariate empirical \chf)
it follows that $G_h$ has covariance function
\beam\label{eq:covG}\lefteqn{
\Gamma((s,t),(s',t'))= \cov(G_h(s,t),G_h(s',t'))}\nonumber \\
&=& \sum_{j\in\mathbb{Z}} \E\big[\big(\ex^{i\,\langle s, X_0\rangle} - \varphi_X(s)\big)\big(\ex^{i\,\langle t,Y_{h}\rangle} -
\varphi_Y(t)\big)\big(\ex^{-i\,\langle s',X_j\rangle} - \varphi_X(-s')\big)\big(\ex^{-i\,\langle t',Y_{j+h}\rangle} - \varphi_Y(-t')\big)\big]\,.
\nonumber\\
\eeam
In the special case when $(X_t)$ and $(Y_t)$ are independent sequences $G_h$ is the same across all $h$
with covariance function
\beao
\Gamma((s,t),(s',t')) = \big(\varphi_X(s-s') - \varphi_X(s)\varphi_X(s')\big)\big(\varphi_Y(t-t')-\varphi_Y(t)\varphi_Y(t')\big)\,.\eeao
Since $G_h$  is centered Gaussian its squared $L^2$-norm $\|G_h\|_\mu^2$ has a weighted $\chi^2$-\ds ; see \cite{kuo:1975}, Chapter 1.
The \ds\ of  $\|G_h\|_\mu^2$ is not tractable and therefore one needs resampling methods for determining its quantiles.
\ere
\bre
Corollary~\ref{cor:adcv} can be extended to the joint \con\ of the \fct\ $n\,T_{n,\mu}^{X,Y}(h)$ at finitely many lags $h$, provided
$X_0$ and $Y_h$ are independent for these lags.
\ere
\bre
Corollary \ref{cor:adcv} does not apply when $X_0$ and $Y_h$ are dependent. Then $n\,T_{n,\mu}^{X,Y}(h)\to \infty$ a.s.
and $n\,R_{n,\mu}^{X,Y}(h)\to \infty$ a.s.
\ere

\section{Auto-distance covariance of fitted residuals from {\rm AR}$(p)$ process}\label{sec:arres}
 An often important problem in time series is to assess the goodness-of-fit of a particular model.  As an illustration, consider a causal autoregressive process of order $p$ ({\rm AR}$(p)$) given by the
difference equations,
\beqo
X_t = \sum_{k=1}^p \phi_k X_{t-k}+Z_t\,,~~~t=0,\pm1,\ldots,
\eeqo
where $(Z_t)$ is an iid sequence with a finite moment $\E[|Z|^\kappa]<\infty$ for some $\kappa>0$.  It is further assumed $Z_t$ has mean 0 if $\kappa\ge1$.
It is often convenient to write the AR($p$) process in the form, 
\beqo
Z_t=X_t-{\vect\phi}^T \bfX_{t-1}\,,
\eeqo
where ${\vect\phi}= (\phi_1,\ldots,\phi_p)^T$, $p\ge 1$  and  $\bfX_t= (X_t,\ldots,X_{t-p+1})^T$. Since the process is assumed causal, we can write
 $X_t=\sum_{j=0}^\infty\psi_j\,Z_{t-j}$ for absolutely summable constants $(\psi_j)$;
see \cite{brockwell:davis:1991}, p.~85.
For convenience, we also write $\psi_j=0$ for $j<0$ and $\psi_0=1$.
\par
The least-squares estimator $\wh {\vect \phi}$ of ${\vect\phi}$ satisfies the relation
\beao
\wh {\vect \phi}- {\vect \phi} &=& \Gamma_{n,p}^{-1}\, \dfrac 1n \sum_{t=p+1}^n \bfX_{t-1} \,Z_t\,,
\eeao
where
\beao
\Gamma_{n,p}&=& \dfrac{1}{n} \sum_{t=p+1}^n \bfX_{t-1}^T \bfX_{t-1}\,.
\eeao
If $\sigma^2={\rm var}(Z_t)<\infty$, we have by the ergodic theorem,
\beam\label{eq:ergodic}
\Gamma_{n,p}\stas \Gamma_p = \big(\gamma_X(i-j)\big)_{1\le i,j\le p}\,,\qquad \mbox{where $\gamma_X(h)= \cov(X_0,X_h)\,,h\in\bbz$}\,.
\eeam
 Causality of the process implies that the partial sum $\sum_{t=p+1}^n \bfX_{t-1} \,Z_t$ is a martingale and applying the martingale \clt\ yields
\beam\label{eq:bfQ}
\sqrt{n} \,\big(\wh {\vect \phi}- {\vect \phi}\big) \std \bfQ\,,
\eeam
where $\bfQ$ is $N({\bf0},\sigma^2 \Gamma_p^{-1})$ distributed.

 The residuals of the fitted model are given by
\beam \label{eq:resid}
\wh Z_t= X_t-  \wh {\vect \phi}^T\,\bfX_{t-1}= \big({\vect \phi}-\wh {\vect \phi} \big)^T\,\bfX_{t-1}+ Z_t\,,\qquad t=p+1,\ldots,n\,.
\eeam
For convenience, we set $\wh Z_t=0$, $t=1,\ldots,p$ since this choice does not influence the \asy\ theory.
Each of the residuals $\wh Z_t$ depends on the estimated parameters and hence the residual process exhibits serial dependence.  Nevertheless, we might expect the test statistic based on the distance covariance function of the residuals given by
\beao
T_{n,\mu}^{\wh Z}(h)= \int_\bbr |C_n^{\wh Z}(s,t)|^2\,\mu(ds,dt)
\eeao
to behave in a similar fashion for the true noise sequence $(Z_t)$.  If the model is a good fit, then we would not expect $T_{n,\mu}^{\wh Z}(h)$ to be extraordinarily large.  As observed by \cite{remillard:2009}, the limit distributions for $T_{n,\mu}^{\wh Z}(h)$ and $T_{n,\mu}^{Z}(h)$ are not the same.  As might be expected, the residuals, which are fitted to the actual data, tend to have smaller distance covariance than the true noise terms for lags less than $p$, if the model is correct.  As a result, one can fashion a goodness-of-fit test based on applying the distance covariance statistics to the residuals.   In the following theorem, we show that the distance covariance based on the residuals has a different limit than for the distance covariance based on the actual noise, if the process has a finite variance.  So in applying a goodness-of-fit test, one must make an adjustment to the limit distribution.  Interestingly,  if the noise has heavy-tails, the limits based on the residuals and the noise terms are the same and no adjustment is necessary.

\bth\label{thm:3}
Consider a causal {\rm AR$(p)$} process with iid noise $(Z_t)$. Assume
\beqq \label{int:res}
\int_{\bbr^2}
\big[(1\wedge |s|^2)\,(1\wedge |t|^2) \,\mu(ds,dt)+  (s^2+t^2)\,\1(|s| \wedge |t|>1)\mu(ds,dt) < \infty.
\eeqq
\begin{enumerate}
\item
If $\sigma^2=\var(Z)<\infty$, then
\beam\label{eq:gxih}
n\,T_{n,\mu}^{\wh Z}(h) \cid \|G_h+ \xi_h\|_\mu^2\,, \qquad\mbox{and}\qquad n\,R_{n,\mu}^{\wh Z}(h) \std \dfrac{\|G_h+\xi_h\|_\mu^2}{
T_\mu^Z(0)}\,,
\eeam
where $(G_h,\xi_h)$ are jointly Gaussian limit random fields on  $\bbr^2$. The covariance structure of $G_h$ is specified
in Remark~\ref{rem:5k} for the \seq\ $((Z_t,Z_{t+h}))$, $\xi_h$ is defined in Lemma~\ref{eq:jj} and the joint limit structure
is described in Lemmas~\ref{lem:1956} and \ref{lem:1957}.
\item
If $Z$ is in the domain of attraction of a stable law of index $\alpha \in (0,2)$, i.e., $\P(|Z|>x) = x^{-\alpha}L(x)$ for $x>0$ and  $L(\cdot)$ is a slowly varying function at $\infty$, and
\beqo
\frac{\P(Z>x)}{\P(|Z|>x)} \to p\, \mbox{\, and \,}\, \frac{\P(Z<-x)}{\P(|Z|>x)} \to 1-p
\eeqo
as $x\to\infty$ for some $p \in [0,1]$ (\cite{feller:1971}, p.~313).
Then we have
\beam\label{eq:gxih2}
n\,T_{n,\mu}^{\wh Z}(h) \cid \|G_h\|_\mu^2\, \qquad\mbox{and}\qquad n\,R_{n,\mu}^{\wh Z}(h) \std \dfrac{\|G_h\|_\mu^2}{
T_\mu^Z(0)}\,,
\eeam
where $G_h$ is a Gaussian limit random field on  $\bbr^2$. The covariance structure of $G_h$ is specified
in Remark~\ref{rem:5k} for the \seq\ $((Z_t,Z_{t+h}))$.
\end{enumerate}
\ethe
The proof is given in Appendix~\ref{sec:thm3}.
\bre
\cite{remillard:2009} mentioned that $T_{n,\mu}^Z(h)$ and $T_{n,\mu}^{\wh Z}(h)$ for an {\rm AR}$(1)$ process have distinct limit processes and he also suggested
the limiting structure in \eqref{eq:gxih}.
\ere
\par

The structure of the limit process in \eqref{eq:gxih} is rather implicit. In applications, one depends on resampling methods.
Relation \eqref{eq:gxih} can be extended to a joint \con\ result for finitely many lags $h$ but the dependence structure of the
limiting vectors is even more involved.  Condition \eqref{int:res}
holds for \pro y \ms s $\mu=\mu_1\times \mu_1$ on $\bbr^2$ with finite second moment but it does not hold for
the  benchmark \ms\ $\mu=\mu_1\times \mu_1$ described in  \eqref{eq:szeklyi}. A reason for this is that $\|\xi_h\|_\mu^2$ is in general
not well defined in this case. 
If  $Z_t$ has characteristic function $\varphi_Z$ then by virtue of 
\eqref{eq:jj},
$\|\xi_h\|_\mu^2$ is finite a.s. if and only if
$$
\int_{-\infty}^\infty |t\varphi_Z(t)|^2 \mu_1(dt)\int_{-\infty}^\infty |\varphi_Z'(s)|^2 \mu_1(ds)<\infty\,.
$$
Now assume that $Z_t$ has a density function $f$ and choose $\mu_1(dt)=c_1t^{-2}dt$. Then by Plancherel's identity, 
the first integral becomes 
$$
\int_{-\infty}^\infty|\varphi_Z(t)|^2 \,dt=c\int_{-\infty}^\infty f^2(t)\,dt\,.
$$
If one chooses $f$ to be a symmetric gamma distribution with shape parameter $\delta\in (0,1/2)$, 
i.e., $f(z)=.5\beta^\delta |z|^{\delta-1}\ex^{-|z|\beta}/\Gamma(\delta)$, then the integral $\int_{-\infty}^\infty f^2(t)dt=\infty$ and hence the limit random variable in \eqref{eq:gxih} cannot be finite.

\medskip
\noindent{\bf AR simulation.}
We illustrate the results of Theorem \ref{thm:3}. First, we generate independent replications of a time series $(X_t)_{t=1,\ldots,1000}$ from a causal {\rm AR}$(10)$ model with $Z_t \sim N(0,1)$ and $$\vect\phi=( -0.140 , 0.038 , 0.304 , 0.078 , 0.069 , 0.013 , 0.019 , 0.039 , 0.148 , -0.062 ).$$
In this and the following examples, we choose the weight measure $\mu$ as the $N(0,0.5)$-\ds\ for which \eqref{int:res} is satisfied.
From the independent replications of the simulated residuals we approximate
the limit \ds\ $\|G_h+\xi_h\|_\mu^2\,/\,T_\mu^Z(0)$ of $n\,R_{n,\mu}^{\wh Z}(h)$
by the corresponding empirical distribution.
\par
\begin{figure} [h]
\begin{center}
\includegraphics[width=17cm]{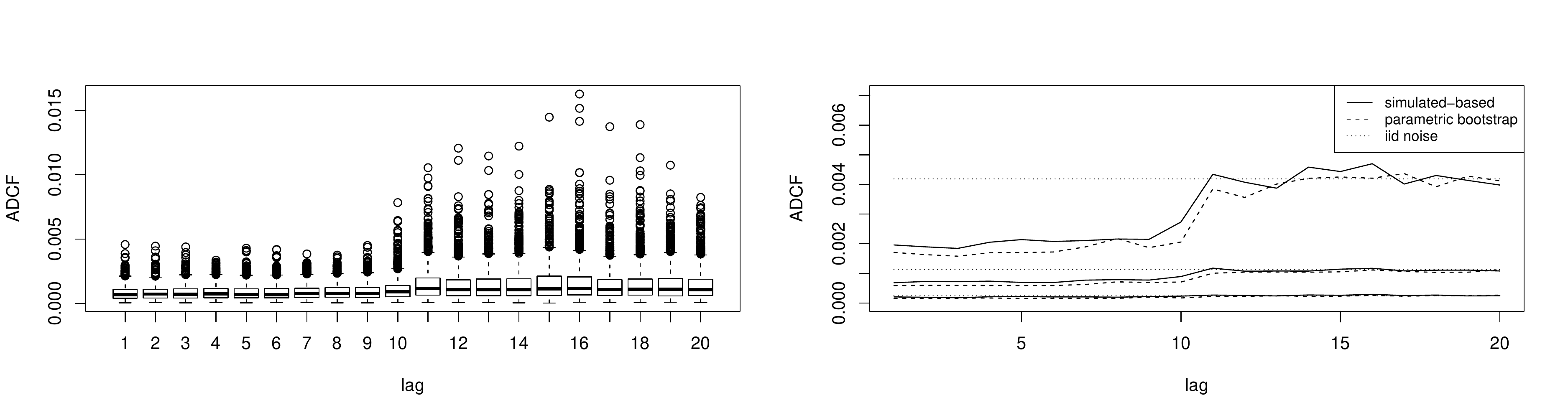}
\end{center}
\caption{Distribution of $n\,R_{n,\mu}^{\wh Z}(h)$, $n=1000$  for the residuals of an AR(10)  process with $N(0,1)$
innovations. Left: Box-plots from 1000 independent replications. Right: $5\%$, $50\%$, $95\%$ empirical quantiles of  $n\,R_{n,\mu}^{\wh Z}(h)$
based on simulated residuals, on resampled residuals and on iid noise, respectively.}
\label{fig:ar}
\end{figure}

The left graph in  Figure \ref{fig:ar} shows the  box-plots for $n\,R_{n,\mu}^{\wh Z}(h)$ based on 1000 replications from the AR(10) model, each with sample size $n=1000$.  As seen from the plots, the distribution at each lag is heavily skewed.
In the right panel of Figure \ref{fig:ar}, we compare the empirical $5\%$, $50\%$, $95\%$
quantiles of $n\,R_{n,\mu}^{\wh Z}(h)$ to those of $n\,R_{n,\mu}^{Z}(h)$,
the scaled ADCF of iid noise, all of which have the same
limit, $\|G_h\|_\mu^2\,/\,T_\mu^Z(0)$.
The asymptotic variance of the ADCF of the residuals is smaller than that of iid noise at initial lags,
and gradually increases at larger lags to the values in the iid case. This behavior is similar to that of the ACF of the residuals of an AR process; see for example Chapter 9.4 of \cite{brockwell:davis:1991}.

Theorem \ref{thm:3} provides a visual tool for testing the goodness-of-fit
of an {\rm AR}$(p)$ model, by examining the serial dependence of the residuals after model fitting.
Under the null hypothesis, we expect $n\,R_{n,\mu}^{\wh Z}(h)$ to be well bounded by the $95\%$ quantiles of the limit distribution $\|G_h+\xi_h\|_\mu^2\,/\,T_\mu^Z(0)$.  For a single time series, this quantity can be approximated using a parametric bootstrap (generating an AR(10) process from the estimated parameters and residuals); see
for example \cite{politis:romano:wolf:1999}.
In the right graph of Figure \ref{fig:ar} we overlay the empirical $5\%$, $50\%$, $95\%$ quantiles of $n\,R_{n,\mu}^{\wh Z}(h)$
estimated from one particular realization of the time series.
As can be seen in the graph, the parametric bootstrap provides a good approximation to the  actual quantiles found via simulation.  On the other hand, the quantiles found by simply bootstrapping the residuals provides a rather poor approximation, at least for the first 10 lags.
\begin{figure} [htbp]
\begin{center}
\includegraphics[width=17cm]{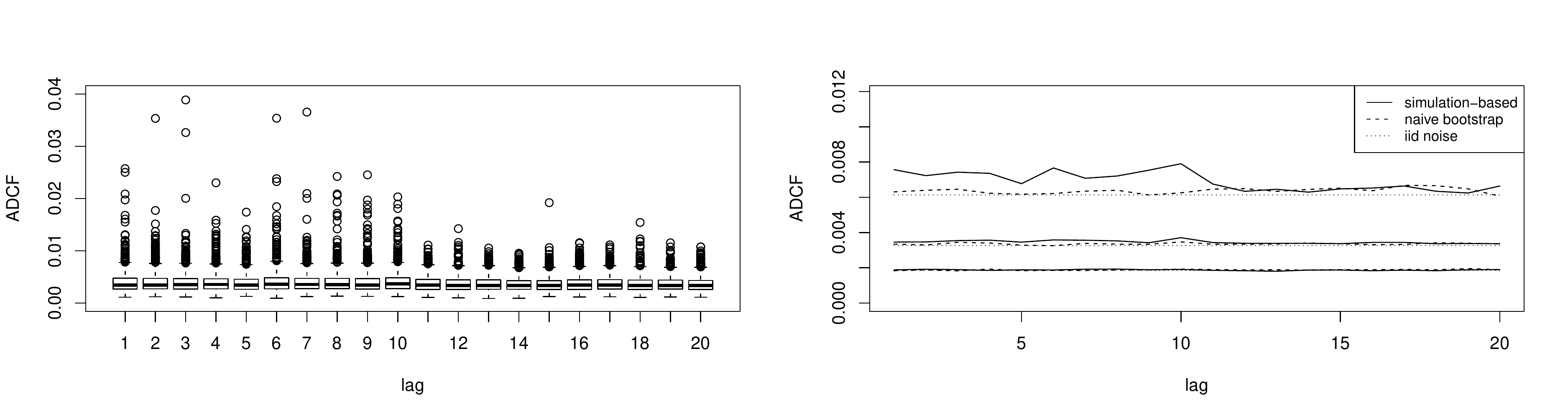}
\end{center}
\caption{Distribution of $n\,R_{n,\mu}^{\wh Z}(h)$ for residuals of AR process with $t_{1.5}$ innovations.
Left: lag-wise box-plots. Right panel: empirical $5\%$, $50\%$, $95\%$ quantiles from simulated residuals, empirical quantiles from resampled residuals, and empirical quantiles from iid noise.}
\label{fig:arht}
\end{figure}

We now consider the same AR(10) model as before, but with noise having a  $t$-distribution with 1.5 degrees of freedom. (Here the noise is in the domain of attraction of a stable distribution with index 1.5.) The left graph of Figure~\ref{fig:arht}
shows the box-plots of $n\ R_{n,\mu}^{\wh Z}(h)$ based on 1000 replications, and the right graph shows the $5\%$, $50\%$, $95\%$ quantiles of
$n\,R_{n,\mu}^{\wh Z}(h)$ and $n\,R_{n,\mu}^{Z}(h)$, both of which have the same limit distribution
$\|G_h\|_\mu^2\,/\,T_\mu^Z(0)$.  
In this case,  the quantiles of $\|G_h\|_\mu^2\,/\,T_\mu^Z(0)$ can be approximated naively by bootstrapping the fitted residuals $({\wh Z_t})$ of the AR model.  The left graph of Figure \ref{fig:arht} overlays the $5\%$, $50\%$, $95\%$
quantiles from bootstrapping with those from the simulations.  The agreement is reasonably good.

\begin{figure} [h]
\begin{center}
\includegraphics[width=17cm]{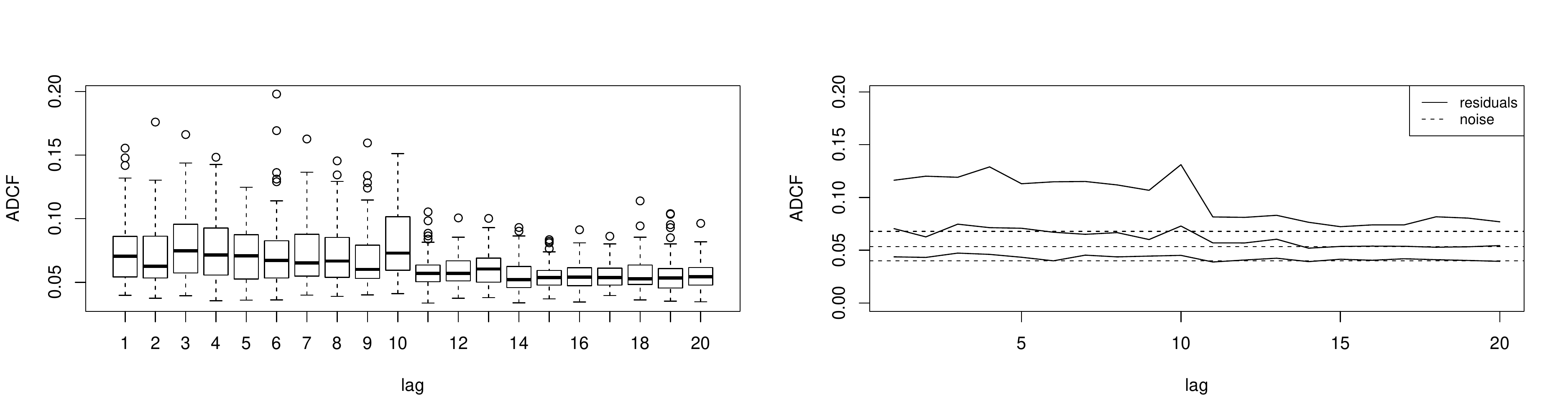}
\end{center}
\caption{Distribution of $n\,R_{n,\mu}^{\wh Z}(h)$, $n=1000$ for residuals of AR process with a symmetric Gamma(0.2,0.5) noise.
Left: box-plots from 500 independent replications. Right panel: empirical $5\%$, $50\%$, $95\%$ quantiles from simulated residuals and from iid noise.}
\label{fig:argamma}
\end{figure}

We next provide an empirical example illustrating the limitation of using the  
measure in \eqref{eq:szeklyi}.  Again, we use the same AR(10) model as before, but with noise now generated from the symmetric gamma distribution with $\delta=.2,\beta=.5$.  The corresponding pair  of graphs with boxplots and quantiles for $n\,R_{n,\mu}^{\wh Z}(h)$ is displayed in Figure~\ref{fig:argamma}.  The $95\%$ quantiles for $n\,R_{n,\mu}^{\wh Z}(h)$ for lags 1-10 are now rather large compared to those of iid noise.

\section{Data Examples}\label{sec:dataex}

\subsection{Amazon daily returns}

\begin{figure} [h]
\begin{center}
\includegraphics[width=17cm]{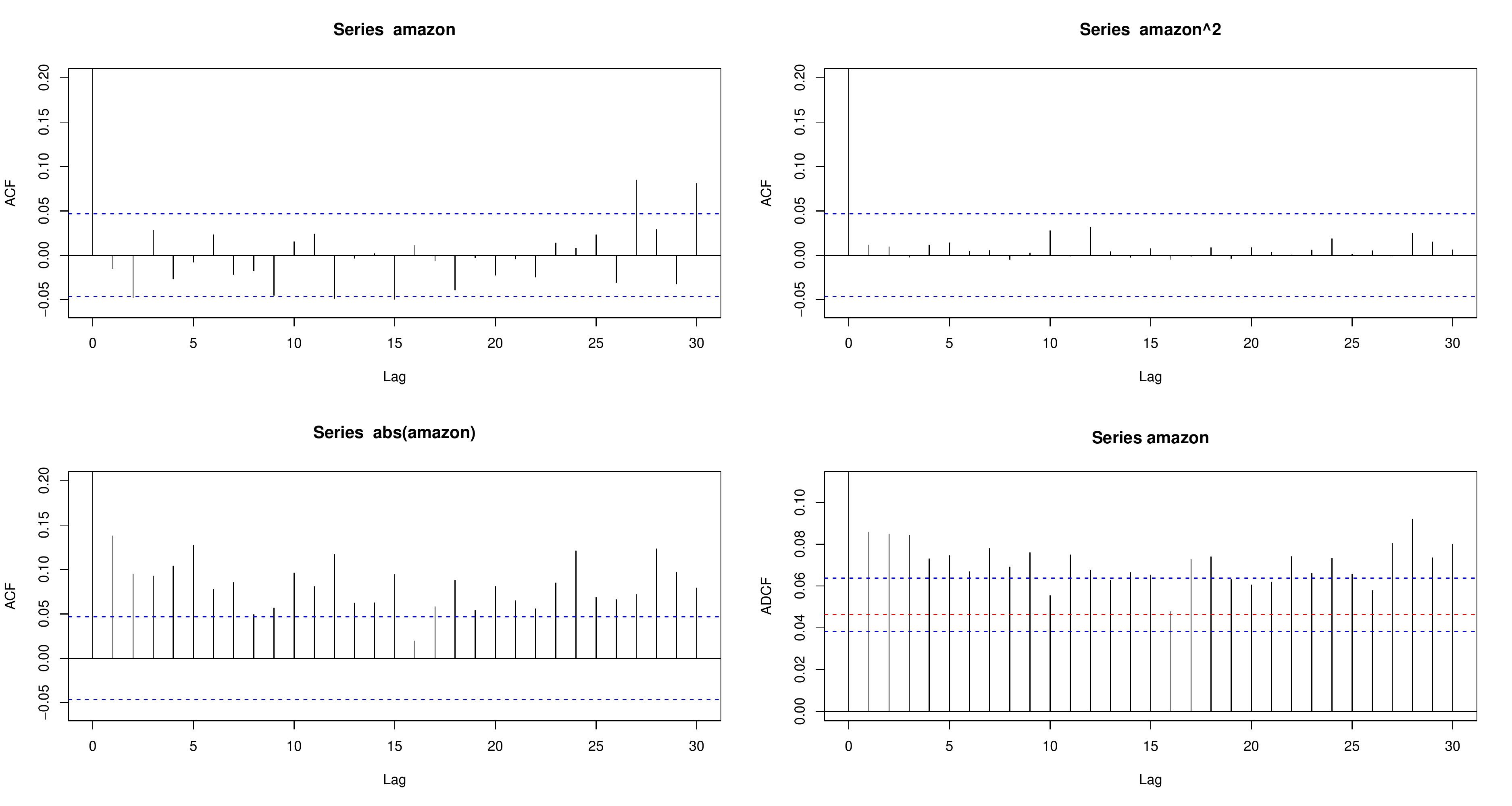}
\end{center}
\caption{ACF and ADCF of daily stock returns of Amazon $(X_t)$  from 05/16/1997 to 06/16/2004. Upper left: ACF of $(X_t)$; Upper right: ACF of $(X^2_t)$; Lower left: ACF of $(|X_t|)$; Lower right: ADCF of $(X_t)$, the $5\%$, $50\%$, $95\%$ confidence bounds of ADCF from randomly permuting the data.}
\label{fig:amazon}
\end{figure}

In this example, we consider the daily stock returns of Amazon from 05/16/1997 to 06/16/2004. 
Denoting the series by $(X_t)$, Figure \ref{fig:amazon} shows the ACF of $(X_t)$, $(X^2_t)$, $(|X_t|)$ and ADCF of $(X_t)$ with weight measure $\mu(ds,dt) = s^{-2}t^{-2}ds dt$. In the right panel, we compare the ADCF with the $5\%$, $50\%$, $95\%$ confidence bounds of the ADCF for iid data, approximated by the corresponding empirical quantiles from 1000 random permutations.  With most financial time series, which are typically uncorrelated, serial dependence can be detected by examining the ACF of the absolute values and squares. Interestingly for the Amazon data,  the ACF of the squared data also fails to pick up any signal. On the other hand, the ADCF has no trouble detecting serial dependence without having to resort to applying any transformation.

\subsection{Wind speed data}

\begin{figure} [h]
\begin{center}
\includegraphics[width=17cm]{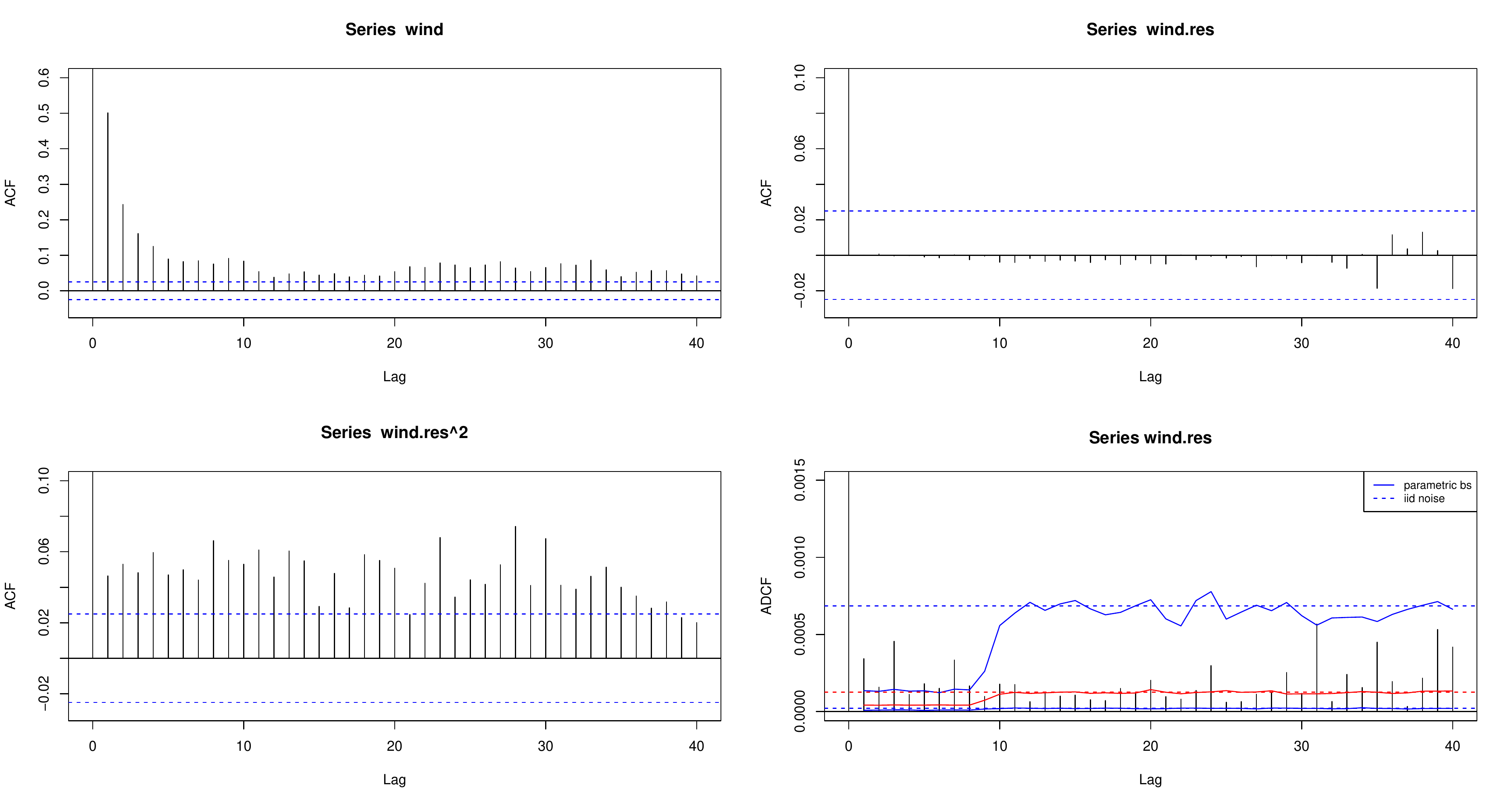}
\end{center}
\caption{ACF and ADCF of Kilkenny wind speed time series and AR(9) fitted residuals. 
Upper left: ACF of the series. 
Upper right: ACF of the residuals. 
Lower left: ACF of the residual squares. 
Lower right: ADCF of the residuals, the $5\%$, $50\%$, $95\%$ confidence bounds of ADCF for fitted residuals from 1000 parametric bootstraps, and that for iid noise from 1000 random permutations. }
\label{fig:kilken}
\end{figure}

For the next example we consider the daily averages of wind speeds at Kilkenny's synoptic meteorological station in Ireland. The time series consists of 6226 observations from 1/1/1961 to 1/17/1978, after which a square root transformation has been applied to stabilize the variance.  This transformation has also been suggested in previous studies (see, for example, \cite{haslett:raftery:1989}). The ACF of the data, displayed in Figure \ref{fig:kilken}, suggests a possible AR model for the data.  An AR(9) model was found to provide the best fit (in terms of minimizing AICC among all AR models) to the data.  The  ACF of the residuals  (see upper right panel in Figure \ref{fig:kilken}) shows that the serial correlation has been successfully removed.  The ACF of the squared residuals and ADCF of the residuals are also plotted in the bottom panels Figure \ref{fig:kilken}.   For computation of the ADCF, we used the N(0,.5) distribution for the weight measure, which satisfies the condition \eqref{int:res}. The ADCF of the residuals is well bounded by the confidence bounds for the ADCF of iid noise, shown by the dotted line in the plot. Without adjusting these bounds for the residuals,  one would be tempted to conclude that the AR model is a good fit. However, the adjusted bounds for the ADCF of residuals, represented by the solid line in the plot and computed using a parametric bootstrap, suggest that some ADCF values among the first 8 lags are in fact larger than expected. Hence this sheds some doubt on the validity of an AR(9) model with iid noise for this data. A similar conclusion can be reached by inspecting the ACF of the squares of the residuals (see lower left panel in Figure \ref{fig:kilken}).

\begin{figure} [h]
\begin{center}
\includegraphics[width=8.5cm]{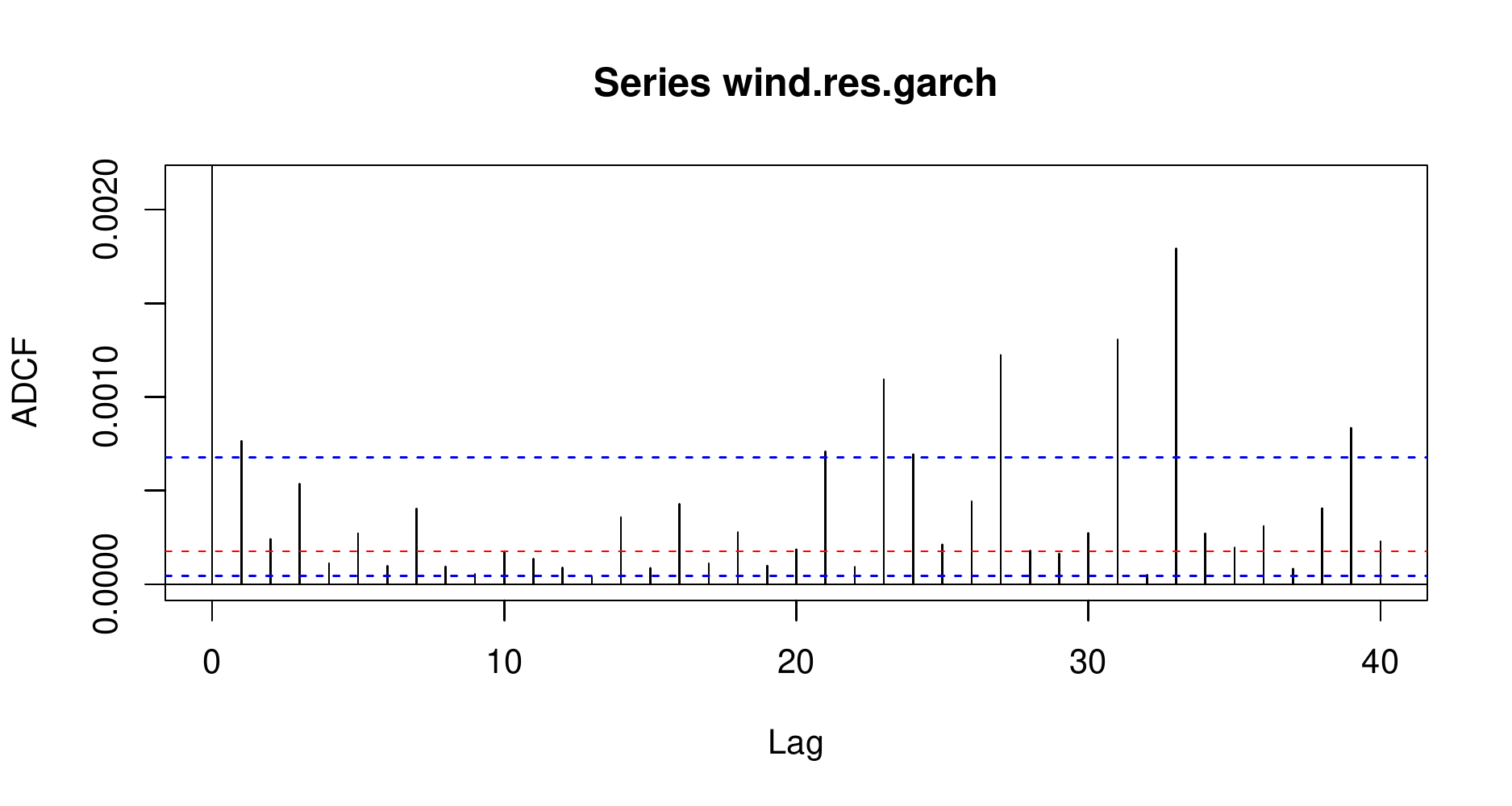}
\end{center}
\caption{ADCF of the residuals of Kilkenny wind speed time series from AR(9)-GARCH fitting and the $5\%$, $50\%$, $95\%$ confidence bounds of ADCF for iid noise from 1000 random permutations.}
\label{fig:kilken2}
\end{figure}

One potential remedy for the lack of fit of the AR(9) model, is to consider a GARCH(1,1) model applied to the residuals.  The GARCH model performs well in devolatilizing the AR-fitted residuals  and no trace of a signal could be detected through the ACF of the GARCH-residuals applied to the squares and absolute values. The ADCF of the devolatilized residuals, seen in Figure \ref{fig:kilken2}, still presents some evidence of dependence. Here the confidence bounds plotted are for iid observations, obtained from 1000 random permutations of the GARCH-residuals and as such do not include an adjustment factor.  Ultimately, a periodic AR model, which allows for periodicity in both the AR parameters and white noise variance might be a more desirable model.

\appendix

\section{Proof of Theorem~\ref{thm:2}}\label{sec:proofoftheorem}
The proof follows from the following lemma.
\ble\label{lem4}
Assume that $\sum_h \alpha_h^{1/r}<\infty$ for some $r>1$ and set $u=2r/(r-1)$.
We also assume  the moment conditions \eqref{eq:moment}  (or \eqref{comd:mom})
for some $\alpha>0$ if $X_0$ and $Y_0$ are independent (dependent).
\begin{enumerate}
\item
For $\alpha\le 2$
there exists a constant $c>0$ \st\ for $\epsilon\in [0,1/2)$,
\beam \label{eq:lem4}
n \,\E[|C_n(s,t)-C(s,t)|^2]\le c\,\big(1\wedge |s|^{\alpha(1+\epsilon)/u} \big)\,\big(1\wedge   |t|^{\alpha(1+\epsilon)/u} \big)\,,
\qquad
n\ge 1\,.
\eeam
\item If $\alpha\in (u/2,u]$ then
$\sqrt{n}(\varphi_{X,Y}^n-\varphi_{X,Y})\std G$ on compact sets $K\subset \bbr^{p+q}$ for some complex-valued mean-zero Gaussian field $G$.
\end{enumerate}
\ele
\bre Notice that $C(s,t)=0$ when $X_0$ and $Y_0$ are independent.
\ere
\begin{proof}
(1) We focus on the proof under the assumption of independence.
At the end, we indicate the changes necessary when $X_0$ and $Y_0$ are dependent.\\
We write
\beao
U_k=\ex^{i\langle s,X_k\rangle}-\varphi_X(s)\,, \qquad V_k=\ex^{i\langle t,Y_k\rangle}-\varphi_Y(t)\,,\qquad k\ge 1\,,
\eeao
where we suppress the dependence of $U_k$ and $V_k$ on $s$ and $t$, respectively.
Then
\beao
n \,\E[|C_n(s,t)|^2]&=& n\,\E\Big|\frac 1 n \sum_{k=1}^n U_k\,V_k-
\dfrac 1n  \sum_{k=1}^n U_k \dfrac 1n \sum_{l=1}^n V_l\Big|^2\\
&\le &2\, n\,\E\Big[\Big|\dfrac 1n \sum_{k=1}^n U_k\,V_k\Big|^2\Big]+ 2\,n\,\E
\Big[\Big|\dfrac 1n  \sum_{k=1}^n U_k \dfrac 1 n\sum_{l=1}^n V_l\Big|^2\Big]\\&=& 2 \,(I_1+I_2)\,.
\eeao
We have by stationarity
\beao
I_1&=& \E [|U_0V_0|^2] + 2 \sum_{h=1}^{n-1} (1-h/n)\,{\rm Re}\, \E [U_0V_0\,\ov{U_hV_h}]\,.
\eeao
Since $U_0$ and $V_0$ are independent $\E[U_0V_0]=0$. In view of the $\alpha$-mixing condition  (see \cite{doukhan:1994}, Section 1.2.2, Theorem 3(a)) we have
\beam\label{eq:ww}
\big|{\rm Re} \E [U_0V_0\,\ov{U_hV_h}]\big|
&\le& c\,\alpha_h^{1/r}\,(\E [|U_0V_0|^u])^{2/u}\\
&=&c\,\alpha_h^{1/r}\,(\E [|U_0|^u])^{2/u} (\E [|V_0|^u])^{2/u}\nonumber\\
&\le &c\,\alpha_h^{1/r}\,(\E [|U_0|^2])^{2/u} (\E [|V_0|^2])^{2/u}\,.\nonumber
\eeam
In the last step we used that $u=2r/(r-1)>2$ and that $\max (|U_0|,|V_0|)\le 2$.
We have for $\alpha\in (0,2]$
\beao
\E[|U_0|^2]
= 1-|\varphi_X(s)|^2\le \E [1 \wedge |\langle s,X-X'\rangle|^\alpha]
\le c\,\big(1\wedge |s|^{\alpha}\big)\,.
\eeao
Therefore and since $\sum_h \alpha_h^{1/r}<\infty$ we have
\beao
I_1&\le &c \, \big(1\wedge |s|^{\alpha}\big)^{2/u} \,
\big(1\wedge |t|^{\alpha}\big)^{2/u}\,.
\eeao
Now we turn to $I_2$.
By the Cauchy-Schwarz inequality and since $|\frac 1 n \sum_{k=1}^n U_k|$ and $|\frac 1 n \sum_{k=1}^n V_k|$ are bounded by 2 we have
\beao
I_2&\le& 2\,n\, \Big(\E \Big|\dfrac 1 n \sum_{k=1}^n U_k\Big|^4\Big)^{1/2}\, \Big(\E \Big|\dfrac 1 n \sum_{k=1}^n V_k\Big|^4\Big)^{1/2}\\
&\le & c\,
\Big(n\,\E \Big|\dfrac 1 n \sum_{k=1}^n U_k\Big|^{2+\delta}\Big)^{1/2}\, \Big(n\,\E \Big|\dfrac 1 n \sum_{k=1}^n V_k\Big|^{2+\delta}\Big)^{1/2}\,,
\eeao
for any $\delta\in [0,2]$. In view of Lemma 18.5.1 in \cite{ibragimov:linnik:1971}
we have for  $\delta\in [0,1)$,
\beao
I_2&\le& c\,
\Big(n\,\E \Big|\dfrac 1 n \sum_{k=1}^n U_k\Big|^{2}\Big)^{(2+\delta)/4}\, \Big(n\,\E \Big|\dfrac 1 n \sum_{k=1}^n V_k\Big|^{2}\Big)^{(2+\delta)/4}\,,
\eeao
Similar arguments as for $I_1$ show that
 \beao
I_2&\le &c \, \big(1\wedge |s|^{\alpha(2+\delta)/4}\big)^{2/u} \,
\big(1\wedge |t|^{\alpha(2+\delta)/4}\,\big)^{2/u}\,.
\eeao
Combining the bounds for $I_1$ and $I_2$, we arrive at \eqref{eq:lem4}.\\[1mm]
Now we indicate the changes necessary when $X_0$ and $Y_0$ are dependent.
We use the notation above and, additionally, write $\wt W_k= U_kV_k-C(s,t)$. We have
\beao
C_n(s,t)-C(s,t) = \frac{1}{n}\sum_{k=1}^n \wt W_k - \frac{1}{n}\sum_{k=1}^n
 U_k \frac{1}{n} \sum_{l=1}^n V_l\,.
\eeao
Then
\beao
n \,\E[|C_n(s,t)-C(s,t)|^2]
&\le &2\, n\,\E\Big[\Big|\dfrac 1n \sum_{k=1}^n \wt W_k\Big|^2\Big]+ 2\,n\,\E
\Big[\Big|\dfrac 1n  \sum_{k=1}^n U_k \dfrac 1 n\sum_{l=1}^n V_l\Big|^2\Big]\\&=& 2 \,(I_1'+I_2)\,.
\eeao
Since $\E[\wt W_0]=0$, we have by stationarity
\beao
I_1'&=& \E [|\wt W_0|^2] + 2 \sum_{h=1}^{n-1} (1-h/n)\,{\rm Re}\, \E [\wt W_0\,\ov {\wt W_h}]\,.
\eeao
Observe that $\E[|\wt W_0|^2] \le 2(\E|U_0|^4\,\E|V_0|^4)^{1/2} + 2|C(s,t)|^2$ and
\beao
 |U_0|^2 &\le& (|\ex^{i\langle s, X_0\rangle}-1|+\E[|1-\ex^{i\langle s,
 X_0\rangle}|])^2 \\
 &\le & c\,(1 \wedge (|s|\,|X_0|)^{\alpha/2})^2 + c\,( 1\wedge
 (|s|^{\alpha/2}\,\E|X_0|^{\alpha/2})  )^2\,.
\eeao
Since $\E [|X_0|^{2\alpha}]<\infty$
 we have $\E[|U_0|^4]\le c\,(1\wedge |s|^{2\alpha})$ and in a similar manner,
 $\E|V_0|^4\le c\,(1\wedge |t|^{2\alpha})$. We also have $|C(s,t)|^2\le c\,(1\wedge
 |s|^\alpha)\,(1\wedge |t|^\alpha)$. Finally, we conclude that
\beao\E [|\wt W_0|^2]
 \le c\,(1\wedge |s|^\alpha)\,(1\wedge |t|^\alpha)\,.
\eeao
With the $\alpha$-mixing condition we obtain
\beao
\big|{\rm Re}\, \E [\wt W_0\,\ov{\wt W_h}]\big|
&\le& c\,\alpha_h^{1/r}\,(\E [|\wt W_0|^u])^{2/u}
\le c\,\alpha_h^{1/r}\,(\E [|\wt W_0|^2])^{2/u}. \nonumber
\eeao
This together with $\sum_h \alpha_h^{1/r}<\infty$ yields
\beao
I_1'&\le &c \, \big(1\wedge |s|^{\alpha}\big)^{2/u} \,
\big(1\wedge |t|^{\alpha}\big)^{2/u}\,.
\eeao
The remaining term $I_2$ can be treated in the same way as in the independent case.
Combining the bounds for $I_1'$ and $I_2$, we arrive at \eqref{eq:lem4}.\\[2mm]
(2) We need an analog of S. Cs\"org\H{o}'s central limit theorem
 (Cs\"org\H{o}\ (1981a,1981b,1981c)) 
empirical \chf\ of an iid multivariate \seq\
with Gaussian limit.
For ease of notation we focus  on the $X$-\seq ; the proof for the $(X,Y)$-\seq\ is analogous and therefore omitted.
The \con\ of the \fidi s of $\sqrt{n}(\varphi_X^n-\varphi_X)$
follows from Theorem~18.5.2 in \cite{ibragimov:linnik:1971}
combined with the Cram\'er-Wold device.
We need to show tightness of the normalized empirical \chf\ on compact sets.
We use the sufficient condition of Theorem 3 in \cite{bickel:wichura:1971}
for multiparameter processes. We evaluate the process on cubes
$(s,t]=\prod_{k=1}^p (s_k,t_k]$, where $s=(s_1\ldots,s_p)$ and
 $t=(t_1,\ldots,t_p)$ and $s_i<t_i$, $i=1,\ldots,p$. The increment of
the normalized empirical
characteristic function on $(s,t]$ is given by
\beam
 I_n(s,t] &=& \sqrt{n}(\varphi_X^n(s,t]-\varphi_X(s,t]) \nonumber\\
 &=& \frac{\sqrt{n}}{n} \sum_{r=1}^n \Big\{
 \sum_{k_{1}=0,1}\cdots\sum_{k_p=0,1}(-1)^{p-\sum_j\,k_j}
 \big(
 \prod_{l=1}^p \ex^{i(s_l+k_l(t_l-s_l))X_{r}^{(l)}} - \E\big[\prod_{l=1}^p \ex^{i(s_l+k_l(t_l-s_l))X_{r}^{(l)}}\big]
\big)
\Big\} \nonumber\\
 &=& \dfrac{1}{\sqrt{n}} \sum_{r=1}^n W_r\,,\label{eq:wr}
\eeam
where $X_r=(X_{r}^{(1)},\ldots,X_{r}^{(p)})$ and
\beao
  W_r=
 \prod_{l=1}^p (\ex^{it_lX_{r}^{(l)}}-\ex^{is_lX_{r}^{(l)}}) - \E \big[\prod_{l=1}^p
 (\ex^{it_lX_{r}^{(l)}}-\ex^{is_l X_{r}^{(l)}})\big]\,.
 \eeao
 We apply the sums $\sum_{k_j=0,1}$ inductively to derive \eqref{eq:wr}.
Observe that
\beao
  \E \big[|I_n(s,t]|^2\big] = \E [|W_0|^2] + 2\, \sum_{h=1}^{n-1} (1-h/n) {\rm Re}\, \E[W_0\,\ov W_h]\,.
\eeao
By the Lipschitz property of trigonometric functions we have for some constant $c>0$ and $\alpha\in (0,2]$,
 \beao
 |\ex^{is_lX_{r}^{(l)}}-\ex^{it_lX_{r}^{(l)}}|^2 &\le & c\, (1\wedge |t_l-s_l|^2
 (X_{r}^{(l)})^2/4) \\
 &\le & c\,(1\wedge |s_l-t_l|^\alpha |X_{r}^{(l)}|^\alpha/4^\alpha)\,.
\eeao
Proceeding as for \eqref{eq:ww} and
noticing that $\alpha\le 2\le u $, we have
\beao
 |\E[W_0\,\ov W_h]|&\le& c\,\alpha_h^{1/r}\,(\E[|W_0|^{u}])^{2/u} \\
 &\le& \alpha_h^{1/r} \prod_{l=1}^p |s_l-t_l|^{2\alpha/u}  \,\big(\E\big[
 \prod_{l=1}^p |X_{0}^{(l)}|^\alpha\big]\big)^{2/u}.
\eeao
Using the summability of $(\alpha_h^{1/r})$ and the moment condition on $X_0$, we may conclude that
\beao
 \E \big[|I_n(s,t]|^2\big] \le \,c\,\prod_{l=1}^p |s_l-t_l|^{2\alpha/u}.
\eeao
If $2 \alpha/u>1$ the condition of Theorem 3 in \cite{bickel:wichura:1971}
yields that the processes $(\sqrt{n} (\varphi_X^n-\varphi_X))$ are tight on compact sets.
\end{proof}

\bpf[Proof of Theorem~\ref{thm:2}(1)] Recall the definition of $K_\delta$ from \eqref{eq:3w}
and that $X_0$ and $Y_0$ are independent.
From Lemma \ref{lem4} and the continuous mapping theorem we have
\beao
\int_{K_\delta} |\sqrt{n}C_n(s,t)|^2\,\mu(ds,dt)\cid \int_{K_\delta} |G(s,t)|^2\,\mu(ds,dt)\,,\qquad\nto \,.
\eeao
From \eqref{eq:j}, \eqref{eq:lem4} and the dominated convergence theorem, for any $\vep>0$, some $\epsilon\in (0,1/2]$ and $\alpha'\le \min(2,\alpha)$,
\beao
&&\lim_{\delta\downarrow 0} \limsup_{\nto}\P \left(\int_{K_\delta^c}|\sqrt{n}\,C_n(s,t)|^2 \,\mu(ds,dt)>\vep\right)\\
&\le&\vep^{-1}\lim_{\delta\downarrow 0} \limsup_{\nto}\int_{K_\delta^c}\E[|\sqrt{n}C_n(s,t)|^2] \,\mu(ds,dt)\,\\
&\le&\lim_{\delta\downarrow 0} \int_{K_\delta^c} c\,\big(1\wedge |s|^{\alpha'(1+\epsilon)/u} \big)\,\big(1\wedge   |t|^{\alpha'(1+\epsilon)/u} \big) \mu(ds,dt)=0\,.
\eeao
\epf
\bpf[Proof of Theorem~\ref{thm:2}(2)] Now we assume that $X_0$ and $Y_0$ are dependent.
We observe that
\beao
 \sqrt{n}\,(T_n(s,t;\mu)-T(s,t;\mu)) = \int_{\R^{p+q}}\sqrt{n}(|C_n(s,t)|^2-|C(s,t)|^2)\,
 \mu(ds,dt).
\eeao
In view of Lemma \ref{lem4}(2) and the a.s. convergence of $C_n$ on compact sets
the continuous mapping theorem implies that for some Gaussian mean-zero process $G'$,
\beao
&&\int_{K_\delta}
 \sqrt{n}\{(C_n(s,t)-C(s,t))\ov{C}_n(s,t)+C(s,t)(\ov{C}_n(s,t)-\ov{C}(s,t))\}\,\mu(ds,dt)\\
 &&\cid
\int_{K_\delta} G'(s,t)\,\mu(ds,dt)\,,\qquad\nto \,,
\eeao
where $G'X(s,t)=2{\rm Re}\{G(s,t)C(s,t)\}$.  We have
\beao
\big||C_n|^2-|C|^2\big|= \big||C_n-C|^2 + 2\,{\rm Re}\, (\ov C\,(C_n-C))\big|\le c\, |C_n-C|\,.
\eeao
By Markov's inequality, \eqref{eq:lem4} and \eqref{eq:j},
\beao
\lefteqn{\lim_{\delta\downarrow 0}
\limsup_{\nto}\P \left(\int_{K_\delta^c}\sqrt{n} \big||C_n(s,t)|^2-|C(s,t)|^2\big|\,\mu(ds,dt)>\vep\right)}\\
&\le&c\,\lim_{\delta\downarrow 0}
\limsup_{\nto}\int_{K_\delta^c} \big(n\,\E[|C_n-C|^2]\big)^{1/2}\,\mu(ds,dt)\,\\
&\le&\lim_{\delta\downarrow 0} \int_{K_\delta^c} c\,\big(1\wedge
|s|^{\alpha'(1+\epsilon)/u} \big)\,\big(1\wedge 
  |t|^{\alpha'(1+\epsilon)/u} \big)\, \mu(ds,dt)=0\,.
\eeao
\epf

\section{Proof of Theorem~\ref{thm:3}}\label{sec:thm3}
We proof the result for the residuals calculated from least square estimates.
One may show that the same result holds for maximum likelihood and Yule-Walker estimates. We start with a joint \clt\ for $C_n^Z$ and $ \wh {\vect \phi}$.
\ble\label{lem:1956}
For every $h\ge 0$,
\beao
\sqrt{n}\, \big(C_n^Z, \wh {\vect \phi}- {\vect \phi}\big)\std (G_h,\bfQ)\,,
\eeao
where the \con\ is in ${\mathcal C}(K)\times \bbr^p$, $K\subset \bbr^2$ is a compact set,
$G_h$ is the limit process of $C_n^Z$ with covariance structure specified in Remark~\ref{rem:5k} for the \seq\ $((Z_t,Z_{t+h}))$,
$\bfQ$ is the limit in \eqref{eq:bfQ},
$(G_h,\bfQ)$ are mean-zero and jointly Gaussian with covariance matrix
\beam\label{eq_Gprocess}
\cov(G_h(s,t),\bfQ)= -\varphi_Z'(s)\,\varphi_Z'(t)\,\Gamma_p^{-1}  \Psi_h\,,\qquad s,t\in\bbr\,,
\eeam
where $\Psi_h=(\psi_{h-j})_{j=1,\ldots,p}$ and $\varphi_Z'$ is the first derivative of $\varphi_Z$.
\ele
\bpf
We observe that, uniformly for $(s,t)\in K$,
\beao
C_n^Z(s,t) &=& \dfrac 1 {n} \sum_{j=1}^{n-h} \ex^{isZ_j+itZ_{j+h}}-
\dfrac 1 {n }\sum_{j=1}^{n-h} \ex^{isZ_j} \dfrac 1n  \sum_{j=1}^{n-h} \ex^{itZ_{j+h}}\\
&=& \dfrac 1 {n} \sum_{j=1}^{n} \big(\ex^{isZ_j}-\varphi_Z(s)\big)\big(\ex^{itZ_{j+h}}-\varphi_Z(t)\big)\\
&&-\dfrac 1 {n }\sum_{j=1}^{n} \big(\ex^{isZ_j}-\varphi_Z(s)\big)\dfrac 1 {n }\sum_{j=1}^{n} \big(\ex^{itZ_j}-\varphi_Z(t)\big)
+ O_\P (n^{-1})\,.
\eeao
In view of the functional \clt\ for the empirical \chf\ of an iid \seq\
(see Cs\"org\H{o}\ (1981a,1981b)) 
we have uniformly for $(s,t)\in K$,
\beao
\sqrt{n}\, C_n^Z(s,t) &=&\dfrac 1 {\sqrt{n}} \sum_{j=1}^{n} \big(\ex^{isZ_j}-\varphi_Z(s)\big)\big(\ex^{itZ_{j+h}}-\varphi_Z(t)\big)
+O_\P(n^{-1/2})\\
&=& I_n(s,t)+O_\P(n^{-1/2})\,.
\eeao
Therefore it suffices to study the \con\ of the \fidi s of $\big(I_n, \sqrt{n}\,(\wh {\vect \phi}- {\vect \phi})\big)$. In view of
\eqref{eq:ergodic} it suffices to show the \con\ of the \fidi s of $\big(I_n, (1/\sqrt{n}) \sum_{j=1}^n \bfX_{j-1}Z_j\big)$. This \con\ follows
by an application of the martingale \clt\ and the Cram\'er-Wold device.
It remains to determine the limiting covariance structure, taking into account the causality of the process $(X_t)$.
We have
\beao
\cov\big(I_n, \dfrac 1 {\sqrt{n}} \sum_{j=1}^n \bfX_{j-1}Z_j\big)&=&\dfrac 1n \E\Big[
\sum_{j=1}^n \sum_{k=1}^n \big(\ex^{isZ_j}-\varphi_Z(s)\big)\big(\ex^{itZ_{j+h}}-\varphi_Z(t)\big) \bfX_{k-1} Z_k\Big]\,.
\eeao
By causality, $X_k$ and $Z_j$ are independent for $k<j$. Hence $\E [(\ex^{isZ_j}-\varphi_Z(s))(\ex^{itZ_{j+h}}-\varphi_Z(t)) X_{l-k} Z_l]$ is non-zero \fif\
$l=j+h$ and $k\le h$, resulting in
\beao
\lefteqn{\E \big[\big(\ex^{isZ_j}-\varphi_Z(s)\big)\big(\ex^{itZ_{j+h}}-\varphi_Z(t)\big)\, X_{l-k}\, Z_l\big]}\\
&=& \E\big[X_{j+h-k} \big(\ex^{isZ_j}-\varphi_Z(s)\big)\big]\,\E \big[Z_{j+h} \big(\ex^{itZ_{j+h}}-\varphi_Z(t)\big)\big]\\
&=& \psi_{h-k} \,\E \big[Z \big(\ex^{is Z}-\varphi_Z(s)\big)\big]\,\E \big[Z \big(\ex^{it Z}-\varphi_Z(t)\big)\big]\\
&=&-\psi_{h-k} \,i\E \big[Z \ex^{is Z}\big]\,i\E \big[Z \ex^{it Z}\big]\\
&=& -\psi_{h-k}\,\varphi_Z'(s)\,\varphi_Z'(t)\,.
\eeao
This implies \eqref{eq_Gprocess}.
\end{proof}
\ble\label{lem:1957}
For every $h\ge 0$,
\beao
\sqrt{n}\, \big(C_n^Z,  C_n^{\wh Z}-C_n^Z\big)\std (G_h, \xi_h)\,,
\eeao
where $(G_h,\bfQ)$ are specified in Lemma~\ref{lem:1956},
\beam\label{eq:jj}
\xi_h(s,t)=t\varphi_Z(t)\,\varphi_Z'(s) \Psi_h^T \bfQ,\qquad (s,t)\in K\,,\eeam the \con\ is in ${\mathcal C}(K,\bbr^2) $, $K\subset \bbr^2$
is a compact set. In particular, we have
\beam\label{eq:gandxi}
\sqrt{n}\, C_n^{\wh Z}\std G_h+\xi_h\,.
\eeam
\ele
\begin{proof}
We observe that, uniformly for $(s,t)\in K$,
\beam
C_n^{\wh Z}(s,t)-C_n^Z(s,t)&=& \dfrac 1n \sum_{j=1}^n \ex^{is Z_j+it Z_{j+h}}\,\big(\ex^{i\,({\vect \phi}- \wh {\vect \phi})^T
(s\bfX_{j-1}+t \bfX_{j+h-1})}-1\big)\nonumber\\
&& +\dfrac 1n \sum_{j=1}^n \big(1- \ex^{i ({\vect \phi}- \wh {\vect \phi})^T s\bfX_{j-1}}\big)\ex^{isZ_j}
\dfrac 1n \sum_{j=1}^n \ex^{itZ_{j+h}}\nonumber\\
&&+ \dfrac 1n \sum_{j=1}^n \ex^{i ({\vect \phi}- \wh {\vect \phi})^Ts \bfX_{j-1}+ isZ_j}\dfrac1n \sum_{j=1}^n \big(1-
\ex^{i ({\vect \phi}- \wh {\vect \phi})^T t \bfX_{j+h-1}}\big)\ex^ {itZ_{j+h}}
+O_\P(n^{-1})\nonumber\\
&=& E_{n1}(s,t)+E_{n2}(s,t)+E_{n3}(s,t)+O_\P(n^{-1})\,. \label{eq:czhat_cz}
\eeam
Write
\beao
\wt  E_{n1}(s,t) &=&  i\,({\vect \phi}- \wh {\vect \phi})^T\,\dfrac 1n \sum_{j=1}^n (s\bfX_{j-1}+
t\, \bfX_{j+h-1})\,\ex^{is Z_j+it Z_{j+h}}\,.
\eeao
In view of the uniform ergodic theorem, \eqref{eq:bfQ} and the causality of $(X_t)$  we have
\beam\label{eq_rt}
\sqrt{n} \wt  E_{n1}(s,t) \std -i\bfQ^T \E \big[(s\bfX_{0}+ t\bfX_{h})\, \ex^{i (s Z_1+ tZ_{h+1})}\big]= -t \varphi_Z(t)\varphi_Z'(s)
\Psi_h^T\bfQ= \xi_h(s,t)\,,
\eeam
where the \con\ is in ${\mathcal C}(K)$.
By virtue of Lemma~\ref{lem:1956} and the mapping theorem we have  the joint \con\ $\sqrt{n}(C_n^Z,\wt E_{n1})\std (G_h,\xi_h)$ in ${\mathcal C}(K,\bbr^2)$.
Denoting the sup-norm in ${\mathcal C}(K)$ by $\|\cdot\|$, it remains
to show that
$
\sqrt{n} \big(\|E_{n2}\| + \|E_{n3}\| + \|E_{n1}-\wt E_{n1}\|\big)\stp 0\,.
$
The proof for $E_{n2}$ and $E_{n3}$ is analogous to \eqref{eq_rt} by observing that the limiting expectation is zero.
We have by a Taylor expansion for some positive constant $c$,
\beao
\sqrt{n} \|E_{n1}(s,t)-\wt E_{n1}(s,t)\|&\le & c\,\big|\sqrt{n} ({\vect \phi}- \wh {\vect \phi})\big|^2
\,\sup_{(s,t)\in K} \dfrac 1 {n^{3/2}}\sum_{j=1}^n\big|s \bfX_{j-1}+t \bfX_{j+h-1}\big|^2\stp 0\,.
\eeao
In the last step we used the uniform ergodic theorem and  \eqref{eq:bfQ}.
\end{proof}

\bpf[Proof of Theorem~\ref{thm:3}(1)]
We proceed as in the proof of Theorem~\ref{thm:2}. By virtue of \eqref{eq:gandxi} and the continuous mapping theorem
we have
\beao
\int_{K_\delta} |\sqrt{n} \,C_n^{\wh Z}(s,t)|^2\,\mu(ds,dt)\cid \int_{K_\delta} |G(s,t)+\xi_h(s,t)|^2\,\mu(ds,dt)\,,\qquad\nto \,.
\eeao
Thus  it remains to show that
\beqq  \label{eq:axes_int}
\lim_{\delta\downarrow 0} \limsup_{\nto}\P\Big(\int_{K_\delta^c}|\sqrt{n}C_n^{\hat Z}(s,t)|^2 \mu(ds,dt)>\vep\Big)= 0\,,\qquad \vep>0\,.
\eeqq
Following the lines of the proof of Theorem~\ref{thm:2}, we have
\beqo
\lim_{\delta\downarrow 0} \limsup_{\nto}\int_{K_\delta^c}\E [|\sqrt{n}C_n^{Z}(s,t)|^2] \,\mu(ds,dt)= 0\,;
\eeqo
see also Remark~\ref{rem:iid}.
Thus it suffices to show
\beqo
\lim_{\delta\downarrow 0} \limsup_{\nto}\P\Big(\int_{K_\delta^c}|\sqrt{n}(C_n^{\hat Z}(s,t)-C_n^{Z}(s,t))|^2 \mu(ds,dt)>\vep\Big)= 0\,,\qquad \vep>0\,.
\eeqo
For convenience we redefine
\[
 C_n^Z = \dfrac 1 {n} \sum_{j=p+1}^{n-h} \ex^{isZ_j+itZ_{j+h}}-
\dfrac 1 {n }\sum_{j=p+1}^{n-h} \ex^{isZ_j} \dfrac 1n  \sum_{j=p+1}^{n-h} \ex^{itZ_{j+h}}\,.
\]
This version does not change previous results for $C_n^Z$.

Using telescoping sums, we have for $\bar{n}=n-p-h$,
\beao
\lefteqn{\frac{\bar{n}}{n} (C_n^{\hat Z}(s,t)-C_n^Z(s,t))}\\ &=& \frac{1}{\bar{n}}\sum_{j=p+1}^{n-h}A_jB_j
-\frac{1}{\bar{n}}\sum_{j=p+1}^{n-h}A_j\frac{1}{\bar{n}}\sum_{j=p+1}^{n-h}B_j
-\frac{1}{\bar{n}}\sum_{j=p+1}^{n-h} U_j \sum_{j=p+1}^{n-h} B_j-
\frac{1}{\bar{n}}\sum_{j=p+1}^{n-h} V_j \sum_{j=p+1}^{n-h} A_j \\
&& +\frac{1}{\bar{n}}\sum_{j=p+1}^{n-h} U_jB_j +
\frac{1}{\bar{n}}\sum_{j=p+1}^{n-h} V_j A_j
=:\sum_{i=1}^6 I_{nj}(s,t),
\eeao
where, suppressing the dependence on $s,t$ in the notation,
\beao
U_j&=&\ex^{isZ_j}-{\varphi}_Z(s)\,,\qquad V_j=\ex^{itZ_{j+h}}-{\varphi}_Z(t)\,,\\
 A_j&=& \ex^{isZ_j}\,(\ex^{is({\vect\phi}-\wh{\vect\phi})'X_{j-1}}-1)\,, \qquad B_j=\ex^{it
 Z_{j+h}}\,(\ex^{is({\vect\phi}-\wh{\vect\phi})'X_{j+h-1}}-1).
\eeao
Write $K_n=|\sqrt{n}({\vect\phi}-\wh {\vect\phi})|$ and $c>0$ for any positive constant
which may differ from line to line.
By Taylor expansions we have
\beao
&&\lefteqn{n\,|I_{n1}(s,t)|^2
\le \Big( \frac{\sqrt{n}}{\bar{n}}\sum_{j=p+1}^{n-h}|A_jB_j|\Big)^2} \\
  &\le& c \,\Big(\frac{\sqrt n}{\bar{n}}\, \sum_{j=p+1}^{n-h} (1\wedge
  |s|\,|{\vect\phi}-\wh{\vect\phi}|\,|X_{j-1}|) \;(1\wedge
  |t|\,|{\vect\phi}-\wh{\vect\phi}||X_{j+h-1}|)\Big)^2 \\
  &\le & c\, \Big( \,
\min\big(\,|s\,t|\, K_n^{2}\,
\frac{1}{\bar{n}^{3/2}}\,\sum_{j=p+1}^{n-h} |X_{j-1}\,X_{j+h-1}|\,,\ |s|
\, K_n\,\frac{1}{\bar{n}}\,\sum_{j=p+1}^{n-h} |X_{j-1}|\,,\ |t| \,
K_n\,\frac{1}{\bar{n}}\,\sum_{j=p+1}^{n-h} |X_{j+h-1}|\big) \Big)^2\,. 
\eeao
The quantities $K_n$ are stochastically bounded. From ergodic theory, $n^{-1}\sum_{j=1}^{n} |X_{j}| = O_\P(1)$ and $n^{-3/2}\sum_{j=1}^{n} |X_{j}\,X_{j+h}| = o_\P(1)$. Hence
\beqo
 n\,|I_{n1}(s,t)|^2 \, \le \, \min(s^2,t^2, (st)^2) \, O_\P(1) \, \le \, \left((1\wedge s^2)\,(1\wedge t^2) + (s^2+t^2) \1( |s|\wedge |t|\ge 1)\right)\, O_\P(1),
\eeqo
where the term $O_\P(1)$ does not depend on $s$ and $t$.
Thus we conclude for $k=1$ that
\beam\label{eq:kkkk}
\lim_{\delta\downarrow 0} \limsup_{\nto}\P\Big(n\,\int_{K_\delta^c}|I_{nk}(s,t)|^2 \,\mu(ds,dt)>\vep\Big)= 0\,,\qquad \vep>0\,.
\eeam
A similar argument yields
\beao
&&n\, |I_{n2}(s,t)|^2 \\
&\le&\left( \frac{\sqrt n}{\bar{n}^2} \sum_{j,k=p+1}^{n-h}|A_j|\, |B_k|\right)^2 \\
 &\le&  \left( \frac{\sqrt n}{\bar{n}^2}\sum_{j,k=p+1}^{n-h} (1\wedge |s|\,|{\vect\phi}-\wh{\vect\phi}|\,|X_{j-1}|)
 (1\wedge |t|\,|{\vect\phi}-\wh{\vect\phi}|\,|X_{k+h-1}|)\right)^2 \\
 &\le& c\, \left(\min\left(|st| \,K_n^2\,  \frac{1}{\bar{n}^{5/2}}\sum_{j,k=p+1}^{n-h}|X_{j-1}\,X_{k+h-1}|\,,\,
|s| \,K_n\, \frac{1}{\bar{n}}\sum_{j=p+1}^{n-h}|X_{j-1}|\,,\,
|t|\,K_n\, \frac{1}{\bar{n}}\sum_{k=p+1}^{n-h}|X_{k+h-1}|\right)\right)^2\\
&\le& \min(s^2,
t^2, (st)^2)\, O_\P(1).
\eeao
Then \eqref{eq:kkkk} holds for $k=2$.
Taylor expansions also yield
\beao
 n\,|I_{n3}(s,t)|^2 &\le&\left( \frac{\sqrt n}{\bar{n}^2} \sum_{j,k=p+1}^{n-h}|U_j|\, |B_k|\right)^2 \\
 &\le& c\, \left( \frac{\sqrt n}{\bar{n}^2}\sum_{j,k=p+1}^{n-h}  (1\wedge \frac{1}{2}|s|\,(|Z_{j}|+\E|Z|))(1\wedge |t|\,|{\vect\phi}-\wh{\vect\phi}|\,|X_{k+h-1}|)
\right)^2 \\
&\le& \min(t^2, (st)^2)\, O_\P(1).
\eeao
This proves \eqref{eq:kkkk} for $k=3$. By a symmetry argument but
with the corresponding bound $\min(s^2, (st)^2)\, O_\P(1)$
, \eqref{eq:kkkk}
for $k=4$ follows as well.
By Taylor expansion, we also have
\beao
n\,|I_{n5}(s,t)|^2 &\le&\left( \frac{\sqrt n}{\bar{n}} \sum_{j=p+1}^{n-h}|U_j|\, |B_j|\right)^2 \\
 &\le& c\, \left( \frac{\sqrt n}{\bar{n}}\sum_{j=p+1}^{n-h}  (1\wedge \frac{1}{2}|s|\,(|Z_{j}|+\E|Z|))(1\wedge |t|\,|{\vect\phi}-\wh{\vect\phi}|\,|X_{j+h-1}|)
\right)^2 \\
&\le&  \min(t^2, (st)^2)\, O_\P(1).
\eeao
We may conclude that \eqref{eq:kkkk} holds for $k=5$. The case $k=6$ follows in a similar way with the corresponding bound
$\min(s^2, (st)^2)\, O_\P(1)$.
\epf
\bpf[Proof of Theorem~\ref{thm:3}(2)]
We follow the proof of Theorem \ref{thm:3}(1) by first showing that
\beqq \label{eq:ht_integrand}
\sqrt{n}\, C_n^{\wh Z}\std G_h
\eeqq
 in ${\mathcal C}(K) $ for $K\subset \bbr^2$ compact, and then \eqref{eq:axes_int}.
The convergence $\sqrt{n}\, C_n^{Z}\std G_h$ in ${\mathcal C}(K) $ continues to hold as in the proof of Theorem \ref{thm:3}
since the conditions in Cs\"org\H{o}\ (1981a,1981b) 
are satisfied if some moment of $Z$ is finite.
For \eqref{eq:ht_integrand} it suffices to show that
\beqq \label{eq:conv_czhat}
\sqrt{n}\, (C_n^{\wh Z} - C_n^Z)\cip 0\,
\eeqq
 in $ {\mathcal C}(K) $.
Recalling the decomposition \eqref{eq:czhat_cz}, we now can show directly that
$\sup_{|s|,|t|\le M}\sqrt{n}|E_{ni}(s,t)| \cip 0$ for any $M>0$ and $i=1,2,3$,
which implies \eqref{eq:conv_czhat}. We focus only on the case $i=1$ to illustrate the method; the cases $i=2,3$ are analogous.
We observe that for $\delta>0$,
\beam\label{eq:lo}
\sup_{|s|,|t|\le M}\sqrt{n}|E_{n1}(s,t)| &\le&
\sup_{|s|,|t|\le M}\sqrt{n}|{\vect \phi}- \wh {\vect \phi}|\,\dfrac 1n \sum_{j=p+1}^{n-h}  |s\bfX_{j-1}+t\, \bfX_{j+h-1}| \nonumber\\
&\le& M\,n^{\frac{1}{\delta}}|{\vect \phi}- \wh {\vect \phi}|\, \ n^{-\frac{1}{\delta} - \frac{1}{2}}\,\sum_{j=1}^{n}
|\bfX_{j}|\,.
\eeam
On the other hand,
under the conditions of Theorem \ref{thm:3}(2)
\cite{hanna:kanter:1977} showed that $\delta>\alpha$,
\beqo
n^{1/\delta}\,(\vect\phi-\wh{\vect\phi}) \stas 0.
\eeqo
For $\alpha\in (1,2)$, $\E [|\bfX|]<\infty$ and since we can choose $\delta=2$ \st\ $1/\delta+1/2=1$.
The ergodic theorem finally yields that the \rhs\ in \eqref{eq:lo} converges to zero a.s.
As regards the case $\alpha\in (0,1]$, we have $\E[|\bfX|^{\alpha-\gamma}]<\infty$ for any small $\gamma$  and
\beqo
\E\Big[\big|n^{-1/\delta - 1/2}\,\sum_{j=1}^{n} |\bfX_{j}|\big|^{\alpha-\gamma}\Big] \le\,
n^{-\,(\alpha-\gamma)(1/\delta+1/2)+1}\,  \E[|\bfX|^{\alpha-\gamma}]\to 0.
\eeqo
If we choose $\delta$ close to $\alpha$ and $\gamma$ close to zero
the \rhs\ in \eqref{eq:lo} converges to zero in \pro y.
\par
Using the same bounds as in part (1), but writing this time $K_n=n^{1/\delta}|{\vect\phi}-\wh{\vect\phi}|$, we have
\beao
n\,|I_{n1}(s,t)|^2
&\le& c
\Big(\min\big(
|s\,t|\,K_n^2\, n^{-1/2-2/\delta}\sum_{j=1}^{n}\,|X_{j-1} X_{j+h-1}|\,,
|s|\,K_n\, n^{-1/\delta-1/2}\sum_{j=0}^{n}\,|X_{j}|\,,\\&&
|t|\,K_n\, n^{-1/\delta-1/2}\sum_{j=0}^{n}\,|X_{j}|
\big)\Big)^2 \\
&\le& c\, \min(\,|s\,t|^2,\,|s|^2 ,\,|t|^2\,)\,\max\big(K_ n^2 \,n^{-1/2-2/\delta}\,\sum_{j=1}^{n}|X_{j-1}X_{j+h-1}|,
K_n\, n^{-1/\delta - 1/2}\,\sum_{j=0}^{n}\,|X_{j}|\Big)^2 \,.
\eeao
The same argument as above shows that $n^{-1/\delta - 1/2}\,\sum_{j=0}^{n}\,|X_{j}|=O_\P(1)$ for $\delta$ close to $\alpha$.
Since $2|X_{j-1} X_{j+h-1}|\le X_{j-1}^2+ X_{j+h-1}^2$ a similar argument shows that $n^{-1/2-2/\delta}\,\sum_{j=1}^{n}|X_{j-1}X_{j+h-1}|=O_\P(1)$.
These facts establish \eqref{eq:kkkk} for $k=1$.
The same arguments show that  bounds analogous to part (1) can be derived for $n\,|I_{nk}(s,t)|^2$ for $k=2,\ldots,6$. We omit further details.
\epf



\end{document}